\documentclass[12pt]{amsart}
\usepackage{amsmath,amsfonts,amstext,amsthm,epsfig}

\newtheorem{thm}{Theorem}[section]
\newtheorem{cor}[thm]{Corollary}
\newtheorem{lem}[thm]{Lemma}
\newtheorem{prop}[thm]{Proposition}
\newtheorem{hyp}[thm]{Hypothesis}
\newtheorem{rem}[thm]{Remark}
\numberwithin{equation}{section}
\newcommand{\eps}{\varepsilon}

\newcommand{\umlaut}{\"}

\DeclareMathOperator{\supp}{supp}

\newcommand{\la}{\left\langle}
\newcommand{\ra}{\right\rangle}
\newcommand{\ral}{\right\rangle_{L^2}}

\newcommand{\Lt}{{L^2}}
\newcommand{\Lto}{{L^2}}
\newcommand{\Dn}{\Delta_\nu}

\DeclareMathOperator{\Lip}{Lip}
\DeclareMathOperator{\Log}{Log}
\DeclareMathOperator{\spec}{spec}

\def\Z{\mathbb Z}
\def\N{\mathbb N}
\def\R{\mathbb R}
\def\Rn{{\mathbb R}^n_x}

\begin{document}

\title[Conditional stability up to the final time]
{Conditional stability up to the final time for backward-parabolic equations\\ with Log-Lipschitz coefficients}%

\author{D. Casagrande, D. Del Santo and M. Prizzi}
\address{D. Casagrande, Dipartimento Politecnico di Ingegneria e Architettura, Universit\`a degli Studi di Udine, Via delle Scienze 206 - 33100 Udine, Italy} \email{daniele.casagrande@uniud.it}
\address{D. Del Santo, Dipartimento di Matematica e Geoscienze, Universit\`a degli Studi di Trieste, Via A. Valerio 12 - 34100 Trieste, Italy} \email{delsanto@units.it}
\address{M. Prizzi, Dipartimento di Matematica e Geoscienze, Universit\`a degli Studi di Trieste, Via A. Valerio 12 - 34100 Trieste, Italy} \email{mprizzi@units.it}
%\thanks{}%

%\date{}%
%\dedicatory{}%
%\commby{}%
% ----------------------------------------------------------------
\begin{abstract}
We prove logarithmic conditional stability up to the final time for backward-parabolic operators whose coefficients are Log-Lipschitz continuous in $t$ and Lipschitz continuous in $x$. The result complements previous achievements of Del Santo and Prizzi (2009) and Del Santo, J\umlaut ah and Prizzi (2015), concerning conditional stability (of a type intermediate between H\umlaut older and logarithmic), arbitrarily closed, but not up to the final time. 
\end{abstract}
\maketitle
% ----------------------------------------------------------------
%%%%%%%%%%%%%%%%%%
\section{Introduction}
%%%%%%%%%%%%%%%%%%%
In real world models, deterministic diffusion processes are often irreversible. Consider for example the heat equation
$$
\partial_t u=\Delta u
$$
with Cauchy data $u(0,x)=u_0(x)$. The forward initial value problem is well posed in an appropriate space of physically meaningful  configurations, but the evolution has a strong regularizing effect, so when one tries to reconstruct an initial configuration  $u(0,x)$ from a final observation $u(T,x)$ at a positive time $T$, one needs to impose regularity conditions on $u(T,x)$, while in general the backward problem with Cauchy data at $T$ has no solution. However, in a physical context an observation at a final time $T$ records the configuration resulting from an {\em actual} evolution, so the problem of existence is less relevant than that of uniqueness and sensitiveness to errors in measurements.
In \cite{John1960} John introduced the notion of {\it well-behaved problem} for ill-posed problems. According to John a problem is {\it
well-behaved} if
{\it ``only a fixed percentage of the significant
digits need be lost in determining the solution from the data"} \cite[p. 552]{John1960}. More precisely,  a problem is
well-behaved if its solutions in a space $\mathcal H$ depend
H\"older continuously on the data belonging to a space $\mathcal
K$, provided the solutions satisfy a prescribed {\it a priori} bound. According to the literature, we call {\it conditional stability} any continuous dependence (possibly weaker than H\"older) which is subordinated to a prescribed {\it a priori} bound.

In this paper we carry on the investigation about conditional stability of backward solutions for a general parabolic equation. For ease of notation we reformulate the problem inverting the sign of the time variable, so we deal with (forward) solutions of the backward-parabolic equation
\begin{equation}
\partial_t u
+\sum_{i,j}\partial_{x_i}(a_{ij}(t,x)\partial_{x_j}u)=0 \label{intro1}
\end{equation}
on the strip $ [0,T]\times {\mathbb R}^n$. We assume throughout
the paper that the matrix $(a_{ij})^n_{i,j=1}$ is symmetric and
positive definite and that the coefficients $a_{ij}$'s are at least Lipschitz continuous in $x$ and H\umlaut older continuous in $t$. These are the standard regularity assumptions which guarantee the (forward) well posedness for forward-parabolic equations in $H^s$, $0\leq s\leq2$ (see e.g. \cite{Amann}).
We denote by $${\mathcal H}:=C^0([0,T], L^2({\mathbb R}^n))\cap C^0([0,T),
H^1({\mathbb R}^n))\cap C^1([0,T), L^2({\mathbb R}^n))$$ the space for admissible solutions of (\ref{intro1}).

In \cite{AN1963} Agmon and Nirenberg proved, among other
things, that the Cauchy problem for (\ref{intro1}) on the interval $[0,T]$ is well-behaved in the space ${\mathcal H}$  with data
in $L^2(\R^n)$ on each subinterval $[0, T']$ with $T'<T$, provided the coefficients $a_{i,j}$'s are
sufficiently smooth with respect to $x$ and Lipschitz continuous
with respect to $t$. In order to achieve their result they developed the so
called {\it logarithmic convexity technique}. The main step
consists in proving that the function $t\mapsto
\log\|u(t,\cdot)\|_{L^2}$ is convex for every solution
$u\in{\mathcal H }$ of (\ref{intro1}). In the same year Glagoleva
\cite{Glagoleva1963} obtained essentially the same result for a concrete
operator like (\ref{intro1}) with time independent coefficients.
Her proof rests on energy estimates obtained through integration
by parts. Some years later Hurd \cite{Hurd1967} developed the technique
of Glagoleva to cover the case of a general equation of type
(\ref{intro1}), with coefficients depending Lipschitz continuously
on time. The results of \cite{AN1963,Glagoleva1963,Hurd1967} can be summarized as
follows:

\medskip

{\bf Theorem A.} {\it Assume the coefficients $a_{ij}$'s are Lipschitz continuous with respect to $t$. For every $T'\in\,\,(0,T)$ and $D>0$ there exist $\rho>0$,
$0<\delta<1$ and $K>0$ such that, if $u\in {\mathcal H}$ is a
solution of (\ref{intro1}) on $[0,T]$ with
$\|u(0,\cdot)\|_{L^2}\leq\rho$ and $\|u(t,\cdot)\|_{L^2}\leq D$ on
$[0,T]$, then $$\sup_{t\in[0,T']}\|u(t,\cdot)\|_{L^2}\leq
K\|u(0,\cdot)\|_{L^2}^\delta.$$  The constants $\rho$, $K$ and
$\delta$ depend only on $T'$ and $D$, on the positivity constant
of  the matrix $(a_{ij})^n_{i,j=1}$, on the $L^\infty$ norms of the coefficients $a_{ij}$'s and of their spatial derivatives, 
and on the Lipschitz constant of the coefficients $a_{ij}$'s with respect to time.}

\medskip

As $T'$ approaches $T$, the constant $K$ above blows up, while $\delta$ decays to $0$, so one cannot expect that solutions are well behaved up to the final time $T$. From the physical point of view, going back to the forward parabolic equation, this means that the reconstruction of the past from observations at the final time $t=T$ worsens more and more as one gets closer to the initial  time $t=0$. 
Yet, as it was proved by various authors (e.g. Imanuvilov and Yamamoto \cite{ImaYama2014}, Yamamoto \cite{Yama2009}, Isakov \cite{Isakov}), some kind of conditional stability for the backward-parabolic equation (\ref{intro1}) {\em up to the final time $T$} can be recovered if one settles for integral estimates rather than pointwise estimates. Moreover, pointwise estimates can be recoverd by imposing stronger {\it a priori} bounds on the solutions. In any case, however, one doesn't get H\umlaut older dependence but only logarithmic dependence on data. The results of \cite{ImaYama2014,Yama2009,Isakov} can be summarized as follows:

\medskip

{\bf Theorem B.} {\it Assume the coefficients $a_{ij}$'s are Lipschitz continuous with respect to $t$. For every $D>0$ there exist $\rho>0$,
$0<\delta\leq1$ and $K>0$ such that, if $u\in {\mathcal H}$ is a
solution of (\ref{intro1}) on $[0,T]$ with
$\|u(0,\cdot)\|_{L^2}\leq\rho$ and $\|u(t,\cdot)\|_{L^2}\leq D$ on
$[0,T]$, then $$\int_0^T\|u(t,\cdot)\|_{L^2}^2\,dt\leq
K\frac1{|\log\|u(0,\cdot)\|_{L^2}|^{2\delta}}.$$  
Moreover, if $\|u(t,\cdot)\|_{H^1}\leq D$ on $[0,T]$, then 
$$\sup_{t\in[0,T]}\|u(t,\cdot)\|_{L^2}\leq
K\frac1{|\log\|u(0,\cdot)\|_{L^2}|^\delta}.$$
The constants $\rho$, $K$ and
$\delta$ depend only on $D$, on the positivity constant
of  the matrix $(a_{ij})^n_{i,j=1}$, on the $L^\infty$ norms of the coefficients $a_{ij}$'s and of their spatial derivatives, 
and on the Lipschitz constant of the coefficients $a_{ij}$'s with respect to time.}

\medskip

In all the above mentioned results, Lipschitz continuity of the coefficients
$a_{ij}$'s with respect to time plays an essential role. The
possibility of replacing Lipschitz continuity by simple continuity
was ruled out by Miller \cite{Miller1973} and more recently by Mandache
\cite{Mandache}. They constructed examples of operators of the form
(\ref{intro1}) which do not enjoy the uniqueness property in
${\mathcal H}$. In the example of Miller the coefficients
$a_{ij}$'s are H\"older continuous in time, while in the more
refined example of Mandache the modulus of continuity $\bar\mu$ of
the coefficients $a_{ij}$'s with respect to time needs only to satisfy
$\int_0^1(1/\bar\mu(s))ds<+\infty$. On the other hand, in
\cite{DSP2005,DSP2015,DSP2021} it was proved that if
$\bar\mu$ satisfies the {\it Osgood condition}, i.e.
$\int_0^1(1/\bar\mu(s))ds=+\infty$, then equation (\ref{intro1}) enjoys
the uniqueness property in ${\mathcal H}$. Therefore it would be
natural to conjecture that if the Osgood condition  is satisfied,
then the Cauchy problem for (\ref{intro1})  is well-behaved in
$\mathcal H$ with data in $L^2(\R^n)$. Unfortunately this is not
true, as shown by a counterexample in \cite{DSP2009}. Nevertheless
if the coefficients
$a_{ij}$'s are Log-Lipschitz continuous in time, it was shown in \cite{DSP2009,DSJP2015} that a weaker conditional stability result holds:

\medskip

{\bf Theorem C.} {\it Assume the coefficients $a_{ij}$'s are Log-Lipschitz continuous with respect to $t$. For every $T'\in\,\,(0,T)$ and $D>0$ there exist $\rho>0$,
$0<\delta<1$ and $K,N>0$ such that, if $u\in {\mathcal H}$ is a
solution of (\ref{intro1}) on $[0,T]$ with
$\|u(0,\cdot)\|_{L^2}\leq\rho$ and $\|u(t,\cdot)\|_{L^2}\leq D$ on
$[0,T]$, then $$\sup_{t\in[0,T']}\|u(t,\cdot)\|_{L^2}\leq K
e^{-N|\log \|u(0,\cdot)\|_{L^2}|^{\delta }}.$$  The constants
$\rho$, $K$, $N$ and $\delta$ depend only on $T'$ and $D$, on the
positivity constant of  the matrix $(a_{ij})^n_{i,j=1}$, on the $L^\infty$ norms of the
coefficients $a_{ij}$'s and of their spatial
derivatives, and on the Log-Lipschitz constant of the coefficients
$a_{ij}$'s with respect to time.}

\medskip

Moreover, in \cite{CDSP2020} a (very feeble) conditional stability result was proved even 
when the coefficients $a_{ij}$'s are just Osgood continuous with respect to $t$, provided they depend only on time.

The proof of Theorem C relies on weighted energy estimates in the spirit
of \cite{Glagoleva1963,Hurd1967,ImaYama2014,Yama2009}, but in order to overcome the obstructions
created by the lack of time differentiability of the coefficients
$a_{ij}$'s it is necessary to introduce a weight function taylored on the modulus of continuity of the $a_{ij}$'s (see Proposition \ref{thm:Energy}), and a microlocal approximation procedure originally developed by Colombini and Lerner in \cite{CL1995} in the context of hyperbolic equations with Log-Lipschitz coefficients.

In this paper we shall exploit the same type of weighted energy estimates to extend Theorem B to the case of parabolic equations whose coefficients are Log-Lipschitz continuous in time (Theorems \ref{L2final} and \ref{Linftyfinal}). Our results can be summarized as follows:

\medskip

{\bf Theorem D.} {\it Assume the coefficients $a_{ij}$'s are Log-Lipschitz continuous with respect to $t$. For every $D>0$ there exist $\rho>0$,
$0<\delta\leq1$ and $K>0$ such that, if $u\in {\mathcal H}$ is a
solution of (\ref{intro1}) on $[0,T]$ with
$\|u(0,\cdot)\|_{L^2}\leq\rho$ and $\|u(t,\cdot)\|_{L^2}\leq D$ on
$[0,T]$, then $$\int_0^T\|u(t,\cdot)\|_{L^2}^2\,dt\leq
K\frac1{|\log\|u(0,\cdot)\|_{L^2}|^{2\delta}}.$$  
Moreover, if $\|u(t,\cdot)\|_{H^1}\leq D$ on $[0,T]$, then 
$$\sup_{t\in[0,T]}\|u(t,\cdot)\|_{L^2}\leq
K\frac1{|\log\|u(0,\cdot)\|_{L^2}|^\delta}.$$
The constants $\rho$, $K$ and
$\delta$ depend only on $D$, on the positivity constant
of  the matrix $(a_{ij})^n_{i,j=1}$, on the $L^\infty$ norms of the coefficients $a_{ij}$'s and of their spatial derivatives, 
and on the Log-Lipschitz constant of the coefficients $a_{ij}$'s with respect to time.}
\medskip

Our results therefore complement the achievements of \cite{DSP2009,DSJP2015}, and {\it en passant} improve them in some crucial technical points related to the regularity of the coefficients $a_{ij}$'s with respect to the $x$ variable (see the discussion in the final part of section 2). Finally, in Section 6 we illustrate some applications of the main results.

%%%%%%%%%%%%%%%%%%%%%%%%%%%%%%%%%%%%%%%%%%%%%%%%%%%%%%%%%%%%%%%%
\section{The weighted energy estimate} \label{sec:WeightedEE}
%%%%%%%%%%%%%%%%%%%%%%%%%%%%%%%%%%%%%%%%%%%%%%%%%%%%%%%%%%%%%%%%
We consider the backward-parabolic equation \begin{equation} \label{BWeq}
 \partial_t u + \sum\limits_{j,k=1}^n \partial_{x_j}(a_{jk}(t,x)\partial_{x_k}u) = 0
\end{equation} on the strip $[0,T] \times \R^n_x$. 
\begin{hyp} \label{mainhypothesis} We assume throughout the paper that:
\begin{itemize}
\item for all $(t,x) \in [0,T] \times \R^n_x$ and for all $j,k = 1, \dots, n$, \begin{equation*}
a_{jk}(t,x) = a_{kj}(t,x);
\end{equation*}
\item there exists $\kappa \in (0,1)$ such that for all $(t,x,\xi) \in [0,T] \times \R^n_x \times \R^n_\xi$, \begin{equation}
\kappa |\xi|^2 \leq \sum\limits_{j,k=1}^n a_{jk}(t,x)\xi_j \xi_k \leq \frac{1}{\kappa} |\xi|^2;
\end{equation}
\item for all $j,k=1,\dots,n$, $a_{jk} \in \Log\Lip([0,T],L^\infty(\R^n_x)) \cap L^\infty([0,T], \Lip(\R^n_x))$. 
\end{itemize}
We set \begin{eqnarray*}
&& A_{LL} := \sup\Big\{ \frac{|a_{jk}(t,x)-a_{jk}(s,x)|}{|t-s|(1+|\log|t-s||)} \mid j,k = 1, \dots, n, \\
&& \qquad \qquad \qquad\qquad\qquad\qquad \qquad t,s \in [0,T],\, x\in  \R^n_x, \, 0 < |s-t| \leq 1 \Big\}, \\[0.3 cm]
&& A := \sup\{ \|\partial_x^\alpha a_{jk}(t,\cdot)\|_{L^\infty} \mid |\alpha| \leq 1, \, t \in [0,T] \}.
\end{eqnarray*}\end{hyp}

\begin{rem}\label{forward} By classical regularity theory for elliptic partial differential equations (see e.g. \cite[Thms. 8.8 and 8.12]{GilTrud}), for each $t\in[0,T]$ the operator $$\mathcal A(t) u:=-\sum\limits_{j,k=1}^n \partial_{x_j}(a_{jk}(t,x)\partial_{x_k}u)$$  is self-adjoint and positive definite in $L^2(\Bbb R^n)$, with domain $H^2(\Bbb R^n)$. Moreover the dependence on $t$ of the operator $\mathcal A(t)$ is better than H\umlaut older continuous, so one can apply the abstract theory of linear parabolic equations (see e.g. \cite[Thm. 4.4.1]{Amann}) and obtain well posedness of the {\em forward} equation $$\partial_t u - \sum\limits_{j,k=1}^n \partial_{x_j}(a_{jk}(t,x)\partial_{x_k}u) = 0$$
in $H^\theta(\Bbb R^n)$ for every $0\leq\theta\leq2$.\end{rem}

For $s>0$, let $\mu(s)=s(1+|\log(s)|)$. For $p \geq 1$, we define \begin{equation*}
\omega(p) := \int_{\frac{1}{p}}^1 \frac{1}{\mu(s)} ds = \log(1+\log p).
\end{equation*} The function $\omega : [1,+\infty) \rightarrow [0,+\infty)$ is bijective and strictly increasing. For $y \in (0,1]$ and $\lambda > 1$, we set $\psi_\lambda(y) = \omega^{-1}(-\lambda \log(y)) = \exp(y^{-\lambda}-1)$ and we define \begin{equation*}
\Phi_\lambda(y) := -\int_y^1 \psi_\lambda(z) dz.
\end{equation*} The function $\Phi_\lambda :(0,1] \rightarrow (-\infty,0]$ is bijective and strictly increasing; moreover, it satisfies \begin{equation} \label{eq:DefPhi}
y\Phi''_{\lambda}(y)= -\lambda(\Phi_\lambda'(y))^2 \mu\big( \frac{1}{\Phi_\lambda'(y)} \big) =
-\lambda\Phi_\lambda'(y)\big( 1 + |\log\big(\frac{1}{\Phi_\lambda'(y)}\big)| \big).
\end{equation} 
In the next lemma, we collect some properties of the functions $\psi_\lambda$ and $\Phi_\lambda$. The proof is left to the reader.

\begin{lem}\label{weight-properties} Let $\zeta>1$. Then, for $y\leq 1/\zeta$, \begin{equation*}
\psi_\lambda(\zeta y) = \exp(\zeta^{-\lambda}-1)(\psi_\lambda(y))^{\zeta^{-\lambda}}.
\end{equation*} Define $\Lambda_\lambda(y) := y \Phi_\lambda(1/y)$. Then the function $\Lambda_\lambda : [1,+\infty) \rightarrow (-\infty,0]$ is bijective and \begin{equation*}
\lim_{z \rightarrow -\infty} -\frac{1}{z}\psi_\lambda\big( \frac{1}{\Lambda_\lambda^{-1}(z)} \big) = +\infty.
\end{equation*}\qed
\end{lem}

We denote by \begin{equation*}
\mathcal H := C^0([0,T], L^2(\R^n_x)) \cap C^0([0,T),H^{1}(\R^n_x)) \cap C^1([0,T), L^2(\R^n_x))
\end{equation*} the space for admissible solutions of \eqref{BWeq}.

\begin{prop}[Weighted energy estimate] \label{thm:Energy} Assume Hypothesis \ref{mainhypothesis} is satisfied. There exists a  constant $\alpha_1>0$ (depending only on $A_{LL}$, $A$ and $\kappa$) and, setting $\alpha := \max\{\alpha_1, T^{-1}\}$, $\sigma := \frac{1}{\alpha}$ and $\tau  := \frac{\sigma}{4}$, there exist constants $\bar{\lambda}>1$, $\bar{\gamma}>0 $ and $M>0$ (depending on $A_{LL}$, $A$, $\kappa$ and $\alpha$, and hence on $T$) such that, for all  $\beta \geq \sigma + \tau$, $\lambda \geq \bar{\lambda}$ and $\gamma \geq \bar{\gamma}$ and whenever $u \in \mathcal H$ is a solution of equation \eqref{BWeq}, the estimate \begin{eqnarray} 
\nonumber && \int_0^s e^{2\gamma t}e^{-2\beta \Phi_\lambda\left( \frac{t+\tau}{\beta} \right)} \|u(t,\cdot)\|_{H^{1-\alpha t}}^2 dt \\
\nonumber && \qquad \qquad \leq M\gamma \Big( (s+\tau)e^{2\gamma s}e^{-2\beta \Phi_\lambda\left( \frac{s+\tau}{\beta} \right)} \|u(s,\cdot)\|_{H^{1-\alpha s}}^2 \\
\label{EnergyEst1} && \qquad \qquad \qquad \qquad + \tau \Phi'_\lambda\left( \frac{\tau}{\beta} \right) e^{-2\beta \Phi_\lambda\left( \frac{\tau}{\beta} \right)} \|u(0,\cdot)\|_{L^2}^2 \Big)
\end{eqnarray} holds for all $s \in [0,\sigma]$.
\end{prop}

\begin{rem} \label{rem:LowerOrder} If one would like to include lower order terms in \eqref{BWeq}, one has to suppose that the corresponding coefficients are $L^\infty$ with respect to $t$ and also $\Lip$ with respect to $x$. The constants in Proposition \ref{thm:Energy}  then will depend also on the norms of the coefficients of the lower order terms. \end{rem}

In \cite{DSP2009} estimate (\ref{EnergyEst1}) was used to deduce the following local conditional stability result:

\begin{thm}[{\cite[Thm.1]{DSP2009}}] \label{thm:DSP2009local} 
Assume Hypothesis \ref{mainhypothesis} is satisfied.
Let $\alpha_1$, $\alpha$ and $\sigma$ be as in Proposition \ref{thm:Energy}. Then there exist constants $\rho$, $\delta$, $K$ and $N$, such that, whenever $u \in \mathcal H$ is a solution of \eqref{BWeq} with $\|u(0,\cdot)\|_{L^2} \leq \rho$, the inequality 
$$\sup\limits_{t \in [0,\sigma/8]} \|u(t,\cdot)\|_{L^2} 
\leq K(1+\|u(\sigma, \cdot)\|_{L^2}) \exp(-N(|\log(\|u(0,\cdot)\|_\Lt)|^\delta)
$$
holds true. The constants $\rho$, $\delta$, $K$ and $N$
depend on $A_{LL}$, $A$, $\kappa$ and $\alpha$, and hence on $T$.\qed
\end{thm}

The fact that $\alpha_1$ is independent of $T$ and  $\sigma=\min\{\alpha_1^{-1}, T\}$ allows one to iterate the local result of Theorem \ref{thm:DSP2009local} a finite number of times, and to obtain conditional stability in the large.

\begin{thm}[{\cite[Thm. 2]{DSP2009}}] \label{thm:DSP2009global}
Assume Hypothesis \ref{mainhypothesis} is satisfied.
Then for all $T' \in (0,T)$ and $D>0$ there exist positive constants $\rho'$, $\delta'$, $K'$ and $N'$, depending only on $A_{LL}$, $A$, $\kappa$, $T$, $T'$ and $D$, such that if $u \in \mathcal H$ is a solution of \eqref{BWeq} satisfying $\sup_{t\in [0,T]}\|u(t,\cdot)\|_{L^2} \leq D$ and $\|u(0,\cdot)\|_{L^2} \leq \rho'$, the inequality \begin{eqnarray*}
\sup\limits_{t \in [0,T']} \|u(t,\cdot)\|_{\Lt} \leq K' \exp\big(-N'|\log(\|u(0,\cdot)\|_{L^2})|^{\delta'}\big)
\end{eqnarray*} holds true. \qed\end{thm}

\begin{rem}Notice that, following Remark \ref{forward}, it would be sufficient to impose an a-priory bound on $\|u(T,\cdot)\|_{L^2}$, which authomatically implies the a-priori bound for $\|u(t,\cdot)\|_{L^2}$, $t\in[0,T]$. \end{rem}

Estimate (\ref{EnergyEst1}) was proved in \cite{DSP2009} when the coefficients $a_{ij}(t,x)$ are of class $C^2$ with respect to $x$ (in this case the constant $A$ contains also the $L^\infty$ norm of the second order spatial derivitaves of the $a_{ij}$'s).  Actually, in \cite{DSP2009} $C^2$ regularity was imposed to overcome a technical difficulty in managing a commutator term appearing in the dyadic decomposition of equation (\ref{BWeq}). However, once estimate (\ref{EnergyEst1}) is achieved, Theorems \ref{thm:DSP2009local} and \ref{thm:DSP2009global} follow directly from it, and the additional regularity in $x$ of the $a_{ij}$'s  plays no role. 

The $C^2$ requirement is somewhat \lq\lq non natural\rq\rq, since Lipschitz continuity in $x$ of the $a_{ij}$'s is sufficient in order that the domain of the operator $-\sum\limits_{j,k=1}^n \partial_{x_j}(a_{jk}(t,x)\partial_{x_k})$ be $H^2(\Bbb R^n)$ (see \cite[ Thms. 8.8 and 8.12]{GilTrud}).

In \cite{DSJP2015} a weaker version of estimate (\ref{EnergyEst1}) was obtained by mean of Bony paraproducts (see \cite{B}), when $C^2$ regularity in $x$ is replaced by the more natural Lipschitz regularity. In this weaker version of (\ref{EnergyEst1}) the spaces $L^2$ and $H^{1-\alpha t}$ were replaced by $H^{-\bar\theta}$ and $H^{1-\bar\theta-\alpha t}$ respectively, where $0<\bar\theta<1$, and the estimate hold for $s\in[0,\frac78\sigma]$, where $\sigma=(1-\bar\theta)/\alpha$ (\cite[Prop. 2.9]{DSJP2015}). Such weaker version of (\ref{EnergyEst1}), together with some nontrivial modifications of the arguments in \cite{DSP2009}, led eventually to recover the continuity results of Theorems \ref{thm:DSP2009local} and \ref{thm:DSP2009global}. However, the weaker weighted energy estimate of \cite{DSJP2015} turns out to be unfit for the pourpose of reaching any kind of stability up to the final time $T$, especially because in that version of the estimate one can not integrate up to $s=\sigma$ in the left hand side of (\ref{EnergyEst1}), but has to stop at $s= \sigma'<\sigma$. Therefore we shall go back to the strong weighted energy estimate (\ref{EnergyEst1}) and demonstrate it in the Lipschitz continuous case, using some ideas contained in \cite{DSJP2015} and performing a more careful and precise analysis of some terms in the paramultiplication procedure.

%%%%%%%%%%%%%%%%%%%%%%%%%%%%%%%%%%%%%%%%%%%%%%%%%%%%%%%%%%%%%%%%
\section{Littlewood-Paley theory and Bony's paraproduct} 
%%%%%%%%%%%%%%%%%%%%%%%%%%%%%%%%%%%%%%%%%%%%%%%%%%%%%%%%%%%%%%%%

In this section, we review some elements of the Littlewood-Paley decomposition which we shall use throughout this paper to define Bony's paraproduct. The proofs which are not contained in this section can be found in \cite{DSP2009}, \cite{DSP2015} and  \cite{Metivier}.

Let $\chi \in C_0^\infty(\R)$ with $0 \leq \chi(s) \leq 1$ be an even function and such that $\chi(s)=1$ for $|s| \leq 11/10$ and $\chi(s)=0$ for $|s| \geq 19/10$. We now define $\chi_k(\xi) = \chi(2^{-k}|\xi|)$ for $k \in \Z$ and $\xi \in \R^n_\xi$. Denoting by $\mathcal F$ the Fourier-transform and by $\mathcal F^{-1}$ its inverse, we define the operators \begin{eqnarray*}
S_{-1}u = 0 &\text{and}&  S_ku = \chi_k(D_x)u = \mathcal F^{-1}(\chi_k(\cdot)\mathcal F(u)(\cdot)), \, k \geq 0, \\
\Delta_0 u = S_0 u &\text{and}& \Delta_k u = S_k u-S_{k-1}u, \, k \geq 1.
\end{eqnarray*} 
We define 
\begin{equation*}
\spec(u) := \supp(\mathcal F(u))%= \overline{\{\xi \in \R^n_\xi : \mathcal F(u)(\xi) \neq 0 \}},
\end{equation*} and we will use the abbreviation $\Delta_k u = u_k$. For $u \in \mathcal S'(\Rn)$, we have \begin{eqnarray*}
u = \lim\limits_{k \rightarrow +\infty} S_k u = \sum\limits_{k \geq 0} \Delta_k u
\end{eqnarray*} in the sense of $\mathcal S'(\R^n_x)$.

We shall make use of the classical 
\begin{prop}[Bernstein's inequalities] \label{PropBernstein} 
Let $u \in \mathcal S'(\R^n_x)$. Then, for $\nu\geq 1$, 
\begin{equation} \label{bern1}
2^{\nu-1}\|u_\nu\|_{\Lt} \leq \|\nabla_x u_\nu\|_{L^2} \leq 2^{\nu+1} \|u_\nu\|_{\Lt}.
\end{equation} The right inequality of \eqref{bern1} holds also for $\nu = 0$.\qed
\end{prop}

In the following two propositions we recall the characterization of the classical Sobolev spaces and Lipschitz-continuous functions via Littlewood-Paley decomposition.

\begin{prop}[{\cite[Lemma 3.2]{CM2008}}] \label{SobolevLP} Let $\theta \in \R$. Then a tempered distribution $u \in \mathcal S'(\Rn)$ belongs to $H^\theta(\Rn)$ iff for all $k \geq 0$, $\Delta_k u \in L^2(\Rn)$
and $$\sum_{k=0}^{+\infty}2^{2k\theta}\|\Delta_k u\|_{L^2}^2<+\infty.$$ Moreover, there exists $C_\theta \geq 1$ such that for all $u \in H^\theta(\R^n_x)$, we have \begin{equation*}
\frac{1}{C_\theta}\|u\|_{H^\theta} \leq\left( \sum_{k=0}^{+\infty}2^{2k\theta}\|\Delta_k u\|_{L^2}^2\right)^{\frac12}\leq C_\theta \|u\|_{H^\theta}.
\end{equation*}
The constant $C_\theta$ remains bounded for $\theta$ in compact subsets of $\Bbb R$.\qed
\end{prop}

\begin{prop}[{\cite[Lemma 3.2]{GR}}] \label{CharLip} A function $a \in L^\infty(\R^n_x)$ belongs to $\Lip(\R^n_x)$ iff \begin{equation*}
\sup\limits_{k \in \N_0} \|\nabla_x (S_k a)\|_{L^\infty} < +\infty.
\end{equation*} Moreover, there exists a positive constant $C$ such that if $a \in \Lip(\R^n_x)$, then \begin{equation*}
\|\Delta_k a\|_{L^\infty} \leq C 2^{-k} \|a\|_{\Lip}, \quad \text{and} \quad \|\nabla_x (S_k a)\|_{L^\infty} \leq C \|a\|_{\Lip}, \end{equation*} where $\|a\|_{\Lip}=\|a\|_{L^\infty}+\|\nabla\|_{L^\infty}$.\qed\end{prop}

Let $a \in L^\infty(\R^n_x)$. Then, Bony's paraproduct of $a$ and $u \in H^\theta(\R^n_x)$ is defined as \begin{equation*}
T_a u = \sum\limits_{k \geq 3} S_{k-3}a \Delta_k u.
\end{equation*} For the proof of our conditional stability result it is essential that $T_a$ is a positive operator. Unfortunately, this is not implied by $a(x) \geq \kappa > 0$. Therefore, we have to modify the paraproduct a little bit. Following \cite[Sect. 3.3.]{CM2008} we introduce the operator \begin{equation} \label{PPDef}
T_a^m u = S_{m-1}aS_{m+2}u + \sum\limits_{k \geq m+3} S_{k-3}a\Delta_k u,
\end{equation} where $m \in \N_0$; note $T^0_a = T_a$. As it will be shown below, the operator $T_a^m$ is a positive operator for positive $a$ provided that $m$ is sufficiently large. 
The next results were proved for $T_a$, but Lemma 3.10 in \cite{CM2008} guarantees that they hold also for $T^m_a$.

\begin{prop} [{\cite[Prop. 5.2.1 and Thms. 5.2.8 and 5.2.9]{Metivier}}] \label{MappingProp}
Let $m\in \N\setminus \{0\}$ and  let $a\in L^\infty(\R^n_x)$. Let $\theta\in \R $. 

Then $T^m_a$ maps $H^\theta$ into $H^\theta$ and there exists $C_{m,\theta}>0$ depending only on $m$ and $\theta$, such that, for all $u\in H^\theta$,
\begin{equation}\label{cont_T}
\| T^m_a u\|_{H^\theta} \leq C_{m,\theta} \|a\|_{L^\infty}\, \|u\|_{H^\theta}.
\end{equation}
The constant $C_{m,\theta}$ can be chosen independent of $\theta$ when $\theta$ belongs to a compact subset of $\Bbb R$.

Let $m\in \N\setminus \{0\}$ and  let $a\in \Lip(\R^n_x)$. Then
\begin{itemize} 
\item $a-T^m_a$ maps $L^2$ into $H^1$ and there exists $C_1>0$ depending only on $m$, such that, for all $u\in L^2$,
\begin{equation}\label{cont_a-T}
\| au -T^m_a u\|_{H^{1}} \leq C_1 \|a\|_{\Lip}\, \|u\|_{L^2};\end{equation}
\item for every $i=1$, \dots, $n$, the mapping $u\mapsto a\partial_{x_i}u-T^m_a\partial_{x_i}u$ extends from $L^2$ to $L^2$, and there exits $C_0>0$ depending only on $m$, such that, for all $u\in L^2$,
\begin{equation}\label{cont_a-T-2}
\| a\partial_{x_i}u-T^m_a\partial_{x_i}u\|_{L^2} \leq C_0 \|a\|_{\Lip}\, \|u\|_{L^2}.\end{equation}
\end{itemize}%\qed
\end{prop}

\begin{cor}\label{interpolation} Let $\theta\in[0,1]$. Then for every $i=1$, \dots, $n$, the mapping $u\mapsto a\partial_{x_i}u-T^m_a\partial_{x_i}u$ extends from $H^\theta$ to $H^\theta$, and  for all $u\in H^\theta$,
\begin{equation}\label{cont_a-T-3}
\| a\partial_{x_i}u-T^m_a\partial_{x_i}u\|_{H^\theta} \leq {C_0}^{1-\theta}{C_1}^{\theta} \|a\|_{\Lip}\, \|u\|_{H^\theta}.\end{equation}
\end{cor}
\begin{proof}
By Proposition \ref{MappingProp} the operator $(a-T^m_a)\partial_{x_j}$ is continuous from $H^0$ to $H^0$ and from $H^1$ to $H^1$. The result follows by interpolation (see e.g. Theorems B.1, B.2 and B.7 in \cite{McLean}).\end{proof}

Next we state a positivity result for $T_a^m$.

\begin{prop}[{\cite[Cor. 3.12]{CM2008}}] \label{prop:Pos} Let $a \in L^\infty(\R^n_x) \cap \Lip(\R^n_x)$ and suppose that $a(x) \geq \kappa >0$ for all $x \in \R^n_x$. Then, there exists a constant $m_0 = m_0(\kappa, \|a\|_{\Lip})$ such that \begin{equation*}
 \la T_a^m u \mid u \ra_{L^2} \geq \frac{\kappa}{2}\|u\|_{L^2}^2,
\end{equation*} for all $u \in L^2(\R^n_x)$ and $m \geq m_0$. A similar result is true for vector-valued functions if $a$ is replaced by a positive symmetric matrix.\qed
\end{prop}

The next proposition is needed since $T_a^m$ is not self-adjoint. However, the operator $(T_a^m-(T_a^m)^\ast) \partial_{x_j}$ is of order $0$ and maps, if $a$ is Lipschitz, $L^2$ continuously into $L^2$.  

\begin{prop}[{\cite[Prop. 3.8 and 3.11]{CM2008}} and  {\cite[Prop. 3.8]{DSP2015}}] \label{Adjoint} Let $m \in \N$, $a \in L^\infty(\R^n_x) \cap \Lip(\R^n_x)$. Then the mapping $u\mapsto (T_a^m-(T_a^m)^\ast) \partial_{x_j}u$ extends from $L^2$ to $L^2$
and there exists a constant $C_{m} > 0$ such that for all $u \in L^2(\R^n_x)$
 \begin{equation*}
\|(T_a^m-(T_a^m)^\ast) \partial_{x_j}u\|_\Lt \leq C_{m}\|a\|_{\Lip}\|u\|_{L^2}.
\end{equation*}\qed \end{prop}

We end this section with a property of the commutators $[\Delta_k,T^m_a]$ which will be crucial in the proof of the weighted energy estimate.

\begin{prop}[{\cite[Prop. 3.7]{DSP2015}}]\label{comm}
Let $m\in \mathbb N\setminus \{0\}$, let $\theta\in\Bbb R$ and  let $a\in \Lip$. Denote by $[\Delta_k, T^m_a]$ the commutator  between $\Delta_k$ and $T^m_a$.

Then there exists $C_{m,\theta}$ depending only on $m$ and $\theta$ such that for all $u\in H^{1-\theta}$,
\begin{equation}\label{comm_T}
(\sum_{k=0}^{+\infty} 2^{-2k\theta}\|\partial_{x_j}([\Delta_k, T^m_a]\partial_{x_k}u)\|^2_{L^2})^{\frac{1}{2}}\leq C_{m,\theta} \|a\|_{\Lip}\|u\|_{H^{1-\theta}}.
\end{equation}
The constant $C_{m,\theta}$ can be chosen independent of $\theta$ when $\theta$ belongs to a compact subset of $\Bbb R$.\qed
\end{prop}

%%%%%%%%%%%%%%%%%%%%%%%%%%%%%%%%%%%%%%%%%%%%%%%%%%%%%%%%%%%%%%%%
\section{Proof of the weighted energy estimate} 
%%%%%%%%%%%%%%%%%%%%%%%%%%%%%%%%%%%%%%%%%%%%%%%%%%%%%%%%%%%%%%%%
For ease of notation, we write the proof only in one space dimension. We divide the proof in several steps.

\medskip
%%%%%%%%%%%%
{ \em - Microlocalization and approximation}
%%%%%%%%%%%%
\smallskip

Let $u\in\mathcal H$ be a solution of (\ref{BWeq})
Let $\alpha\geq T^{-1}$, $\sigma=1/\alpha$, $\tau=\sigma/4$, $\gamma>0$, $\lambda>1$, $\beta\geq\sigma+\tau$. For $t\in[0,\sigma]$
define
$w(t,x) = e^{\gamma t}e^{-\beta\Phi_\lambda\left( \frac{t+\tau}{\beta} \right)}u(t,x)$. Then $w$ satisfies 
\begin{equation*}
\partial_t  w-\gamma w + \Phi_\lambda'\left( \frac{t+\tau}{\beta} \right)w + \partial_x(a(t,x)\partial_x w) = 0.
\end{equation*}   
Now we add and subtract $\partial_x T_a^m \partial_x w$, where $T_a^m$  is the paramultiplication operator defined in \eqref{PPDef}, with $m\geq m_0(\kappa, A)$, according to the positivity result of Proposition \ref{prop:Pos}. We obtain \begin{equation} \label{eq:P1}
\partial_t w -\gamma w + \Phi_\lambda'\left( \frac{t+\tau}{\beta} \right)w + \partial_x(T_a^m\partial_x w) + \partial_x((a-T_a^m)\partial_x w) = 0.
\end{equation} We set $u_\nu = \Delta_\nu u$, $w_\nu = \Delta_\nu w$ and $v_\nu = 2^{-\alpha t\nu}w_\nu$. Then the function $v_\nu$ satisfies \begin{equation} \label{eq:P2} \begin{aligned}
& \partial_t v_\nu = \gamma v_\nu - \Phi_\lambda'\left( \frac{t+\tau}{\beta} \right)v_\nu - \partial_x(T_a^m\partial_x v_\nu) - \alpha\log2 \nu v_\nu \\
& \qquad - 2^{-\alpha t\nu} \partial_x([\Dn, T_a^m]\partial_x w) - 2^{-\alpha t\nu} \Delta_\nu \partial_x((a-T_a^m)\partial_x w).
\end{aligned}\end{equation} 
Now we make the scalar product of \eqref{eq:P2} with $(t+\tau)\partial_t v_\nu$ in $L^2(\Bbb R_x)$ and obtain \begin{equation} \label{eq:Mod} \begin{aligned}
& (t+\tau)\|\partial_t v_\nu(t)\|_\Lto^2 = \gamma(t+\tau) \la v_\nu \mid \partial_t v_\nu(t)\ral \\
& \qquad \quad -(t+\tau)\la \Phi_\lambda'\left( \frac{t+\tau}{\beta}\right)v_\nu(t) \mid \partial_t v_\nu(t) \ral \\
& \qquad \quad -(t+\tau)\la \partial_x(T_a^m\partial_x v_\nu(t)) \mid \partial_t v_\nu(t) \ral \\
& \qquad \quad -\alpha\log2(t+\tau)\nu \la v_\nu(t) \mid \partial_t v_\nu(t) \ral \\
& \qquad \quad -(t+\tau) 2^{-\alpha t\nu} \la \partial_x([\Dn, T_a^m]\partial_x w(t)) \mid \partial_t v_\nu(t) \ral \\
& \qquad \quad -(t+\tau)2^{-\alpha t\nu} \la \Delta_\nu \partial_x((a-T_a^m)\partial_x w(t)) \mid \partial_t v_\nu(t) \ral.
\end{aligned}\end{equation}

To proceed further, we need to regularize the coefficient $a(t,x)$ with respect to $t$. We take a regular mollifier, i.e. an even, non-negative $\rho \in C_0^\infty(\R)$ with $\supp(\rho) \subseteq [-\frac{1}{2},\frac{1}{2}]$ and $\int_\R \rho(s) ds = 1$. For $\eps \in (0,1]$, we set \begin{equation*}
a_\eps(t,x) = \frac{1}{\eps} \int_\R a(s,x)\rho\left( \frac{t-s}{\eps} \right) ds.
\end{equation*} A straightforward computation shows that for all $\eps \in (0,1]$, we have \begin{eqnarray}
&& a_\eps(t,x) \geq \kappa > 0 \\
&& |a_\eps(t,x)-a(t,x)| \leq A_{LL} \eps(|\log\eps|+1)
\end{eqnarray} as well as \begin{equation*}
|\partial_t a_\eps(t,x)| \leq A_{LL}\|\rho'\|_{L^1(\R)}(|\log\eps|+1)
\end{equation*} for all $(t,x) \in [0,T] \times \R_x$. From these properties of $a_\eps(t,x)$ and by Proposition \ref{MappingProp}, we immediately get \begin{lem} \label{lem:Approx} Let $m \in \N_0$ and $u \in L^2(\R^n_x)$. Then \begin{equation*}
 \|(T_a^m-T_{a_\eps}^m)u\|_\Lto \leq C_{m,0} A_{LL}\eps(|\log\eps|+1) \|u\|_\Lto \end{equation*}
 and
 \begin{equation*} \|T_{\partial_t a_\eps}^m u\|_\Lto \leq C_{m,0} A_{LL} \|\rho'\|_{L^1(\R)}(|\log\eps|+1)\|u\|_{\Lt}.
\end{equation*} \qed\end{lem}

We set \begin{equation*}
a_\nu(t,x) := a_\eps(t,x), \, \text{ with } \, \eps = 2^{-2\nu}.
\end{equation*}

We replace $T_a^m$ by $T_{a_\nu}^m + T_a^m-T_{a_\nu}^m$ in the third term of the right hand side of \eqref{eq:Mod} and we obtain \begin{equation} \label{eq:Mod1} \begin{aligned}
& (t+\tau)\|\partial_t v_\nu(t)\|_\Lto^2 = \gamma(t+\tau) \la v_\nu(t) \mid \partial_t v_\nu(t)\ral \\
& \qquad \quad -(t+\tau)\la \Phi_\lambda'\left( \frac{t+\tau}{\beta}\right)v_\nu(t) \mid \partial_t v_\nu(t) \ral \\
& \qquad \quad -(t+\tau)\la \partial_x(T_{a_\nu}^m\partial_x v_\nu(t)) \mid \partial_t v_\nu(t) \ral  \\
& \qquad \quad -(t+\tau)\la \partial_x((T_a^m-T_{a_\nu}^m)\partial_x v_\nu(t)) \mid \partial_t v_\nu(t) \ral \\
& \qquad \quad -\alpha\log2(t+\tau) \nu \la v_\nu(t) \mid \partial_t v_\nu(t) \ral \\
& \qquad \quad -(t+\tau) 2^{-\alpha t\nu} \la \partial_x([\Dn, T_a^m]\partial_x w(t)) \mid \partial_t v_\nu(t) \ral \\
& \qquad \quad -(t+\tau)2^{-\alpha t\nu} \la \Delta_\nu \partial_x((a-T_a^m)\partial_x w(t)) \mid \partial_t v_\nu(t) \ral.
\end{aligned}\end{equation}

Now we replace $\partial_t v_\nu(t)$ in the term \begin{equation*}
-\alpha\log2(t+\tau) \nu \la v_\nu(t) \mid \partial_t v_\nu(t) \ral
\end{equation*} by the expression on the right hand side of \eqref{eq:P2} and we obtain 

\begin{equation} \label{eq:Mod2} \begin{aligned}
& -\alpha \log2(t+\tau) \nu \la v_\nu(t) \mid \partial_t v_\nu(t) \ral = \\
&\quad -\alpha \gamma \log2(t+\tau)\nu \|v_\nu(t)\|_\Lto^2 \\
&\quad +\alpha\log2(t+\tau)\Phi_\lambda'\left( \frac{t+\tau}{\beta} \right)\nu \|v_\nu(t)\|_\Lto^2 \\
&\quad +\alpha\log2(t+\tau)\nu \la v_\nu(t) \mid \partial_x T_a^m\partial_x v_\nu(t) \ral \\
&\quad +\alpha^2(\log2)^2 (t+\tau) \nu^2 \|v_\nu(t)\|_\Lto^2 \\
&\quad +\alpha \log2(t+\tau)\nu 2^{-\alpha t\nu}  \la v_\nu(t) \mid \partial_x\big( [\Dn,T_a^m]\partial_x w(t) \big) \ral \\
&\quad +\alpha \log2(t+\tau)\nu 2^{-\alpha t\nu}  \la v_\nu(t) \mid \Dn \partial_x\big( (a-T_a^m)\partial_x w(t) \big) \ral.
\end{aligned}\end{equation}

By \eqref{eq:Mod1} and \eqref{eq:Mod2}, we obtain \begin{equation*}\begin{aligned}
&(t+\tau)\|\partial_t v_\nu(t)\|_\Lto^2 = \\
&\qquad \quad \gamma(t+\tau) \la v_\nu(t) \mid \partial_t v_\nu(t) \ral \\
&\qquad \quad -(t+\tau) \Phi_\lambda'\left( \frac{t+\tau}{\beta}\right) \la v_\nu(t) \mid \partial_t v_\nu(t) \ral \\
&\qquad \quad -(t+\tau)\la \partial_x(T_{a_\nu}^m\partial_x v_\nu(t)) \mid \partial_t v_\nu(t) \ral \\
&\qquad \quad -(t+\tau)\la \partial_x((T_a^m-T_{a_\nu}^m)\partial_x v_\nu(t)) \mid \partial_t v_\nu(t) \ral \\
&\qquad \quad +\alpha\log2(t+\tau)\Phi_\lambda'\left( \frac{t+\tau}{\beta} \right)\nu \|v_\nu(t)\|_\Lto^2 \\
&\qquad \quad +\alpha\log2(t+\tau)\nu \la v_\nu(t) \mid \partial_x T_a^m\partial_x v_\nu(t) \ral \\
&\qquad \quad +\alpha^2(\log2)^2 (t+\tau) \nu^2 \|v_\nu(t)\|_\Lto^2 \\
&\qquad \quad -\alpha\gamma \log2(t+\tau) \nu \|v_\nu(t)\|_\Lto^2 \\
&\qquad \quad +\alpha \log2(t+\tau)\nu 2^{-\alpha t\nu} \la v_\nu(t) \mid \partial_x\big( [\Dn,T_a^m]\partial_x w(t) \big) \ral \\
&\qquad \quad +\alpha \log2(t+\tau)\nu 2^{-\alpha t\nu} \la v_\nu(t) \mid \Dn \partial_x\big( (a-T_a^m)\partial_x w(t) \big) \ral \\
&\qquad \quad -(t+\tau) 2^{-\alpha t\nu} \la \partial_x([\Dn, T_a^m]\partial_x w(t)) \mid \partial_t v_\nu(t) \ral \\
&\qquad \quad -(t+\tau) 2^{-\alpha t\nu} \la \Delta_\nu \partial_x((a-T_a^m)\partial_x w(t)) \mid \partial_t v_\nu(t) \ral. \end{aligned} \end{equation*}

A straightforward computation using Leibnitz derivation rule with respect to $t$ yields \begin{equation*}
\gamma(t+\tau)\la v_\nu(t) \mid \partial_t v_\nu(t) \ral = \frac{\gamma}{2}\frac{d}{dt}\Big( (t+\tau)\|v_\nu(t)\|_\Lto^2 \Big) - \frac{\gamma}{2}\|v_\nu(t)\|_\Lto^2
\end{equation*} 
and  
\begin{equation*} \begin{aligned}
&-(t+\tau)\Phi_\lambda'\left( \frac{t+\tau}{\beta} \right) \la v_\nu(t) \mid \partial_t v_\nu(t) \ral = \\
&\qquad -\frac{1}{2}\frac{d}{dt}\Big( (t+\tau)\Phi_\lambda'\left( \frac{t+\tau}{\beta} \right) \|v_\nu(t)\|_\Lto^2 \Big) + \frac{1}{2}\frac{t+\tau}{\beta}\Phi_\lambda''\left( \frac{t+\tau}{\beta} \right)\|v_\nu(t)\|^2_\Lto \\
&\qquad + \frac{1}{2}\Phi_\lambda'\left( \frac{t+\tau}{\beta} \right)\|v_\nu(t)\|^2_{L^2(\R^n)}.
\end{aligned}\end{equation*} 

Next we consider the term $-(t+\tau)\la \partial_x(T_{a_\nu}^m\partial_x v_\nu(t)) \mid \partial_t v_\nu(t)\ral$. From \eqref{PPDef} it can be seen that $\partial_t T_{a_\nu}^m = T_{\partial_t a_\nu}^m + T_{a_\nu}^m\partial_t$. A simple computation then shows that\begin{equation*} \begin{aligned}
&-(t+\tau)\la \partial_x(T_{a_\nu}^m\partial_x v_\nu(t)) \mid \partial_t v_\nu(t)\ral = \\
& \qquad  \frac{1}{2}\frac{d}{dt}\Big( (t+\tau)\la T_{a_\nu}^m \partial_x v_\nu(t) \mid \partial_x v_\nu(t) \ral \Big) \\
& \qquad - \frac{1}{2}\la  T_{a_\nu}^m \partial_x v_\nu(t) \mid \partial_x v_\nu(t) \ral \\
& \qquad -\frac{1}{2}(t+\tau)\la T_{\partial_t a_\nu}^m \partial_x v_\nu(t) \mid \partial_x v_\nu(t) \ral \\
& \qquad -\frac{1}{2}(t+\tau)\la \partial_t \partial_x v_\nu(t) \mid ((T_{a_\nu}^m)^\ast-T_{a_\nu}^m)\partial_x v_\nu(t) \ral.
\end{aligned} \end{equation*} 

Eventually we obtain the identity \begin{equation} \label{eq:EnEstimate} \begin{aligned}
&(t+\tau)\|\partial_t v_\nu(t)\|_{L^2}^2 = \\
& \qquad \frac{\gamma}{2}\frac{d}{dt}\left( (t+\tau)\|v_\nu(t)\|_\Lto^2 \right) - \frac{\gamma}{2}\|v_\nu(t)\|_\Lto^2 \\
& \qquad -\frac{1}{2}\frac{d}{dt}\left( (t+\tau)\Phi_\lambda'\left( \frac{t+\tau}{\beta} \right)\|v_\nu(t)\|_\Lto^2 \right) \\
& \qquad +\frac{1}{2}\Phi_\lambda'\left( \frac{t+\tau}{\beta} \right)\|v_\nu(t)\|_\Lto^2 + \frac{1}{2}\frac{t+\tau}{\beta}\Phi_\lambda''\left( \frac{t+\tau}{\beta} \right)\|v_\nu(t)\|_\Lto^2 \\
& \qquad - (t+\tau)\la \partial_x((T_a^m-T_{a_\nu}^m)\partial_x v_\nu(t)) \mid \partial_t v_\nu(t) \ral \\
& \qquad + \frac{1}{2}\frac{d}{dt}\left( (t+\tau)\la T_{a_\nu}^m \partial_x v_\nu(t) \mid \partial_x v_\nu(t) \ral \right) \\
& \qquad - \frac{1}{2}\la T_{a_\nu}^m \partial_x v_\nu(t) \mid \partial_x v_\nu(t) \ral \\
& \qquad - \frac{1}{2}(t+\tau)\la T_{\partial_t a_\nu}^m\partial_x v_\nu(t) \mid  \partial_x v_\nu(t) \ral \\
& \qquad - \frac{1}{2}(t+\tau)\la \partial_t \partial_x v_\nu(t) \mid ((T_{a_\nu}^m)^\ast - T_{a_\nu}^m)\partial_x v_\nu(t) \ral \\
& \qquad - \alpha\gamma \log2 (t+\tau) \nu \|v_\nu(t)\|_\Lto^2 \\
& \qquad + \alpha\log2(t+\tau)\Phi_\lambda'\left( \frac{t+\tau}{\beta} \right)\nu \|v_\nu(t)\|^2_\Lto \\
& \qquad - \alpha\log2(t+\tau)\nu \la \partial_x v_\nu(t) \mid T_a^m \partial_x v_\nu(t) \ral \\
& \qquad + \alpha^2(\log2)^2 (t+\tau)\nu^2 \|v_\nu(t)\|^2_{L^2} \\
& \qquad + \alpha\log2(t+\tau)\nu 2^{-\alpha t\nu}\la v_\nu(t) \mid \mathcal X_\nu(t) \ral \\
& \qquad - (t+\tau)2^{-\alpha t\nu}\la \mathcal X_\nu(t) \mid \partial_t v_\nu(t) \ral,
\end{aligned}\end{equation} where we have set \begin{eqnarray*}
\mathcal X_\nu(t) &:=&  \left( \partial_x([\Dn,T_a^m]\partial_x w(t)) + \Dn(\partial_x((a-T_a^m)\partial_x w(t))) \right).
\end{eqnarray*}

\medskip
%%%%%%%%%%
{\em - Estimates for $\nu = 0$}
%%%%%%%%%%
\smallskip

In what follows, we denote by $C^{(1)}$, $C^{(2)}$, $C^{(3)}$, \dots \, positive constants which depend only on $A_{LL}$, $A$ and $\kappa$.

Setting $\nu = 0$, we get from \eqref{eq:EnEstimate} \begin{equation*} \begin{aligned}
&(t+\tau)\|\partial_t v_0(t )\|_{L^2}^2 = \\
& \qquad \frac{\gamma}{2}\frac{d}{dt}\left( (t+\tau)\|v_0(t )\|_\Lto^2 \right) - \frac{\gamma}{2}\|v_0(t )\|_\Lto^2 \\
& \qquad -\frac{1}{2}\frac{d}{dt}\left( (t+\tau)\Phi_\lambda'\left( \frac{t+\tau}{\beta} \right)\|v_0(t )\|_\Lto^2 \right) \\
& \qquad +\frac{1}{2}\Phi_\lambda'\left( \frac{t+\tau}{\beta} \right)\|v_0(t )\|_\Lto^2 + \frac{1}{2}\frac{t+\tau}{\beta}\Phi_\lambda''\left( \frac{t+\tau}{\beta} \right)\|v_0(t )\|_\Lto^2 \\
& \qquad - (t+\tau)\la \partial_x((T_a^m-T_{a_0}^m)\partial_x v_0(t )) \mid \partial_t v_0(t ) \ral \\
& \qquad + \frac{1}{2}\frac{d}{dt}\left( (t+\tau)\la T_{a_0}^m \partial_x v_0(t ) \mid  \partial_x v_0(t ) \ral \right) \\
& \qquad - \frac{1}{2}\la T_{a_0}^m \partial_x v_0(t ) \mid \partial_x v_0(t ) \ral - \frac{1}{2}(t+\tau)\la \partial_x v_0(t ) | T_{\partial_t a_0}^m \partial_x v_0(t ) \ral \\
& \qquad - \frac{1}{2}(t+\tau)\la \partial_t \partial_x v_0(t ) \mid ((T_{a_0}^m)^\ast - T_{a_0}^m)\partial_x v_0(t )  \ral \\
& \qquad - (t+\tau) \la  \mathcal X_0(t )  \mid \partial_t v_0(t ) \ral. \\
\end{aligned}\end{equation*} 

By Proposition \ref{prop:Pos} we have
\begin{equation*}
- \frac{1}{2}\la T_{a_0}^m \partial_x v_0(t ) \mid \partial_x v_0(t )\ral \leq -\frac\kappa8\|v_0(t)\|_{L^2}^2.
\end{equation*}
Using Propositions \ref{PropBernstein}, \ref{MappingProp} and Lemma \ref{lem:Approx}, for $N_1$, $N_2>0$, we get 
\begin{equation*}
 |\la \partial_x v_0(t ) \mid T^m_{\partial_t a_0} \partial_x v_0(t ) \ral| \leq C^{(1)} \|v_0\|^2_\Lt, 
\end{equation*}

\begin{equation*}|\la T_{a-a_0}^m\partial_x v_0(t ) \mid \partial_x\partial_t v_0(t ) \ral| \leq C^{(2)} N_1 \|v_0(t )\|_\Lt^2 + \frac{1}{N_1} \|\partial_t v_0(t )\|_\Lt^2, 
\end{equation*}
and
\begin{equation*}
 |\la ((T_{a_0}^m)^\ast - T_{a_0}^m)\partial_x v_0(t )\mid \partial_t \partial_x v_0(t ) \ral| \leq C^{(3)}N_2 \|v_0(t )\|_\Lt^2 + \frac{1}{N_2} \|\partial_t v_0(t )\|_\Lt^2.
\end{equation*} 
Now, we choose $N_1$ and $N_2$ so large that 
\begin{equation*}
\frac{1}{N_1} + \frac{1}{N_2} - \frac{1}{2} < 0
\end{equation*} 
and $\bar{\gamma}$ so large that 
\begin{equation*}
-\frac{\gamma}{4} + \left(C^{(1)} + C^{(2)} N_1 + C^{(3)} N_2\right)(\sigma + \tau) < 0
\end{equation*} 
for $\gamma \geq \bar{\gamma}$. With this choice, the term 
\begin{equation*} 
C^{(1)} (t+\tau) \|v_0(t )\|_{\Lt}^2 + C^{(2)} N_1(t+\tau)\|v_0(t )\|_\Lt^2+ C^{(3)}N_2(t+\tau)\|v_0(t )\|_\Lt^2
\end{equation*} 
is absorbed by $-\frac{\gamma}{4} \|v_0(t )\|_\Lt^2$, and the term 
\begin{equation*}
\frac{1}{N_1} (t+\tau) \|\partial_t v_0(t )\|_\Lt^2 + \frac{1}{N_2} (t+\tau) \|\partial_t v_0(t )\|_\Lt^2
\end{equation*} 
is absorbed by $-\frac{1}{2}(t+\tau)\|\partial_t v_0(t )\|_\Lt^2$. Hence, we get 
\begin{equation*} 
\begin{aligned} 
& \frac12(t+\tau)\|\partial_t v_0(t )\|_\Lt^2 \\
& \qquad \leq \frac{\gamma}{2} \frac{d}{dt}\left( (t+\tau)\|v_0(t )\|_\Lt^2 \right) - \frac{\gamma}{4}\|v_0(t )\|_\Lt^2
 + \frac{1}{2} \Phi_\lambda'\left( \frac{t+\tau}{\beta} \right) \|v_0(t )\|_\Lt^2\\ 
 &\qquad-\frac\kappa8\|v_0(t)\|_{L^2}^2 
-\frac{1}{2}\frac{d}{dt}\left( (t+\tau) \Phi'_\lambda\left( \frac{t+\tau}{\beta} \right)\|v_0(t )\|_\Lt^2 \right) \\
 &\qquad+ \frac{1}{2} \frac{t+\tau}{\beta}\Phi_\lambda''\left( \frac{t+\tau}{\beta} \right) \|v_0(t )\|_\Lt^2 \\
& \qquad + \frac{1}{2} \frac{d}{dt} \left( (t+\tau)\la T_{a_0}^m \partial_x v_0(t ) \mid \partial_x v_0(t ) \ra_\Lt \right) - (t+\tau) \la  \mathcal X_0  \mid \partial_t v_0(t ) \ral.
\end{aligned} 
\end{equation*} 
Further, we recall that $\Phi$ satisfies equation \eqref{eq:DefPhi}, i.e. 
\begin{equation*}
y\Phi''_{\lambda}(y)= -\lambda(\Phi_\lambda'(y))^2 \mu\Big( \frac{1}{\Phi_\lambda'(y)} \Big) =
-\lambda\Phi_\lambda'(y)\left( 1 + |\log\Big(\frac{1}{\Phi_\lambda'(y)}\Big)| \right)
\end{equation*} 
for $\lambda > 1$. From this, we see that \begin{equation*}
\frac{1}{2} \Phi_\lambda'\left( \frac{t+\tau}{\beta} \right) \|v_0(t )\|_\Lt^2 + \frac{1}{2} \frac{t+\tau}{\beta}\Phi_\lambda''\left( \frac{t+\tau}{\beta} \right) \|v_0(t )\|_\Lt^2 < 0,
\end{equation*} 
and thus we get 
\begin{equation*} \begin{aligned}
&\frac{\gamma}{8}\|v_0(t )\|_\Lt^2 \\
&\qquad\leq -\frac{1}{2}(t+\tau)\|\partial_t v_0(t )\|_\Lt^2 + \frac{\gamma}{2} \frac{d}{dt}\left( (t+\tau)\|v_0(t )\|_\Lt^2 \right) - \frac{\gamma}{8}\|v_0(t )\|_\Lt^2\\
 &\qquad-\frac\kappa8\|v_0(t)\|_{L^2}^2
+ \frac{1}{2} \frac{d}{dt} \left( (t+\tau)\la T_{a_0}^m \partial_x v_0(t ) \mid  \partial_x v_0 \ra_\Lt \right) \\
&\qquad- (t+\tau) \la  \mathcal X_0  \mid \partial_t v_0(t ) \ral -\frac{1}{2}\frac{d}{dt}\left( (t+\tau) \Phi'_\lambda\left( \frac{t+\tau}{\beta} \right)\|v_0(t )\|_\Lt^2 \right).
\end{aligned}\end{equation*} 
Integrating in $t$ over $[0,s] \subseteq [0,\sigma]$, we obtain \begin{equation*} \begin{aligned}
&\frac{\gamma}{8}\int_0^s \|v_0(t )\|_\Lt^2 dt \\
& \qquad \leq (\frac{\gamma}{2} + C^{(4)})(s+\tau)\|v_0(s )\|_\Lt^2 + \frac{1}{2} \tau \Phi_\lambda'\left( \frac{\tau}{\beta} \right) \|v_0(0 )\|_\Lt^2 \\
& \qquad \quad - \frac{\gamma}{8} \int_0^s \|v_0(t )\|_\Lt^2 dt -\frac\kappa8\int_0^s\|v_0(t)\|_{L^2}^2dt\\
&\qquad\quad-\frac{1}{2} \int_0^s (t+\tau)\|\partial_t v_0(t )\|_{\Lt}^2 dt 
 -\int_0^s (t+\tau)\la  \mathcal X_0(t ) \mid \partial_t v_0(t ) \ra_\Lt dt,
\end{aligned}\end{equation*} 
where we have used the estimates 
\begin{eqnarray*}
|\la \partial_x v_0(s )| T_{a_0}^m\partial_x v_0(s ) \ral| \leq 2 C^{(4)}\|v_0(s)\|_\Lt^2
\end{eqnarray*} 
and 
\begin{equation*}
\la \partial_x v_0(0 ) | T_{a_0}^m \partial_x v_0(0) \ral \geq \frac{\kappa}{2}\|\partial_x v_0(0 )\|_\Lt^2,
\end{equation*} 
which follow from propositions  \ref{MappingProp} and \ref{prop:Pos} respectively.

\medskip
%%%%%%%%%%%
{\em - Estimates for $\nu \geq 1$}
%%%%%%%%%%%
\smallskip

Now, we consider \eqref{eq:EnEstimate} for $\nu \geq 1$. From Lemma \ref{lem:Approx} and Proposition \ref{Adjoint}, for $N_3$ and $N_4>0$, we obtain 
\begin{equation} 
\label{auxest4} 
\begin{aligned}
&|\la (T_a^m-T_{a_\nu}^m)\partial_x v_\nu(t) \mid  \partial_x\partial_t v_\nu(t) \ral| \\
& \qquad \qquad \leq C_{a,m}^{(5)} N_3 \nu^2  \|v_\nu(t)\|_\Lt^2 + \frac{1}{N_3}\|\partial_t v_\nu(t)\|_\Lt^2\\
& \qquad \qquad \leq C_{a,m}^{(5)} N_3 \nu 2^{2\nu} \|v_\nu(t)\|_\Lt^2 + \frac{1}{N_3}\|\partial_t v_\nu(t)\|_\Lt^2
\end{aligned} 
\end{equation} 
and 
\begin{equation} \label{auxest3} \begin{aligned}
|\la \partial_x v_\nu(t) \mid T_{\partial_t a_\nu}^m \partial_x v_\nu(t)\ral| \leq C^{(6)}_{a,m}\nu 2^{2\nu}\|v_\nu(t)\|^2_\Lt,
\end{aligned} \end{equation} 
as well as 
\begin{equation} \label{auxest2} \begin{aligned}
& |\la ((T_{a_\nu}^m)^\ast-T_{a_\nu}^m)\partial_x v_\nu(t) \mid \partial_t \partial_x v_\nu(t) \ral| \\
& \qquad \qquad \leq C^{(7)}_{a,m}N_4 2^{2\nu} \|v_\nu(t)\|_\Lt^2 + \frac{1}{N_4}\|\partial_t v_\nu(t)\|_\Lt^2
\end{aligned} \end{equation}  
Using again the positivity estimate in Proposition \ref{prop:Pos} as well as Proposition \ref{PropBernstein}, 
we obtain 
\begin{equation} \label{auxest1} \begin{aligned}
& - \alpha\log2(t+\tau)\nu \la \partial_x v_\nu(t) \mid T_{a}^m \partial_x v_\nu(t) \ral \\
& \qquad \qquad \qquad \qquad \leq -\alpha {\frac{\kappa\log2}4}(t+\tau)\nu 2^{2\nu} \| v_\nu(t)\|^2_{\Lt}.
\end{aligned} \end{equation} 
Now, we choose $N_3$ and $N_4$ so large that \begin{equation*}
\frac{1}{N_3} + \frac{1}{N_4} - \frac{1}{2} < 0,
\end{equation*}
and $\alpha_1$ large enough such that \begin{equation*}
- \frac{\alpha_1}{4} {\frac{\kappa\log2}4} + N_3C_{a,m}^{(5)} + C_{a,m}^{(6)} + C_{a,m}^{(7)}N_4 < 0,
\end{equation*} and we set $\alpha := \max\{T^{-1}, \alpha_1\}$. 
 With this choice, we get \begin{equation} \label{finest1} 
 \begin{aligned}
& \frac{\gamma}{4}\|v_\nu(t)\|_\Lto^2 + \frac{1}{2} (t+\tau)\|\partial_t v_\nu(t)\|_{L^2}^2 \\
& \qquad \leq \frac{\gamma}{2}\frac{d}{dt}\left( (t+\tau)\|v_\nu(t)\|_\Lto^2 \right) - \frac{\gamma}{4}\|v_\nu(t)\|_\Lto^2 \\
& \qquad\quad -\frac{1}{2}\frac{d}{dt}\left( (t+\tau)\Phi_\lambda'\left( \frac{t+\tau}{\beta} \right)\|v_\nu(t)\|_\Lto^2 \right) \\
& \qquad\quad +\frac{1}{2}\Phi_\lambda'\left( \frac{t+\tau}{\beta} \right)\|v_\nu(t)\|_\Lto^2 + \frac{1}{2}\frac{t+\tau}{\beta}\Phi_\lambda''\left( \frac{t+\tau}{\beta} \right)\|v_\nu(t)\|_\Lto^2 \\
& \qquad\quad + \frac{1}{2}\frac{d}{dt}\left( (t+\tau)\la  T_{a_\nu}^m \partial_x v_\nu(t) \mid \partial_x v_\nu(t) \ral \right) \\
& \qquad\quad - \alpha\gamma \log2 (t+\tau) \nu \|v_\nu(t)\|_\Lto^2 -\frac{1}{2}\la T_{a_\nu}^m  \partial_x v_\nu(t) \mid \partial_x v_\nu(t) \ral \\
& \qquad\quad + \alpha\log2(t+\tau)\Phi_\lambda'\left( \frac{t+\tau}{\beta} \right)\nu \|v_\nu(t)\|^2_\Lto \\
& \qquad\quad + \alpha^2(\log2)^2 \nu^2 (t+\tau)\|v_\nu(t)\|_{L^2} -\frac{3\alpha}{4}  {\frac{\kappa\log2}4} (t+\tau) \nu 2^{2\nu} \|v_\nu(t)\|_\Lt^2 \\
& \qquad\quad + \alpha\log2\nu 2^{-\alpha t\nu}(t+\tau)\la v_\nu(t) \mid \mathcal X_\nu(t) \ral \\
& \qquad\quad - (t+\tau)2^{-\alpha t\nu} \la \mathcal X_\nu(t) \mid \partial_t v_\nu(t) \ral.
\end{aligned}\end{equation} 
Since $y\Phi_\lambda''(y)=-\lambda \Phi_\lambda'(y)(1+|\log(\Phi_\lambda'(y))|)$, if we take $\lambda\geq\bar\lambda > 2$, we have \begin{equation*}
\frac{1}{4} \frac{t+\tau}{\beta} \Phi''_\lambda\left( \frac{t+\tau}{\beta} \right) \leq - \frac{1}{2}\Phi_\lambda'\left( \frac{t+\tau}{\beta} \right),
\end{equation*} 
and hence, the term $\frac{1}{2}\Phi_\lambda'\left( \frac{t+\tau}{\beta} \right)\|v_\nu(t)\|_\Lto^2$ in \eqref{finest1} is absorbed by the term $\frac{1}{4}\frac{t+\tau}{\beta}\Phi_\lambda''\left( \frac{t+\tau}{\beta} \right)\|v_\nu(t)\|_\Lto^2$. Now we need to absorb \begin{equation} \label{term1}
\alpha\log2(t+\tau)\Phi_\lambda'\left( \frac{t+\tau}{\beta} \right)\nu \|v_\nu(t)\|^2_\Lto.
\end{equation} 
There are two terms in \eqref{finest1} that will help to achieve this. One is \begin{equation} \label{term2}
-\frac{\alpha}{4}  {\frac{\kappa\log2}4} (t+\tau) \nu 2^{2\nu} \|v_\nu(t)\|_\Lt^2
\end{equation} 
and the other one is \begin{equation} \label{term3}
\frac{1}{4}\frac{t+\tau}{\beta}\Phi_\lambda''\left( \frac{t+\tau}{\beta} \right)\|v_\nu(t)\|_\Lto^2.
\end{equation} 
Let $\kappa' = \min\{4\log2,  {\frac{\kappa\log2}4}\}$. If $\nu \geq {1\over 2\log 2}\log\left(\frac{4\log2}{\kappa'}\Phi_\lambda'\Big( \frac{t+\tau}{\beta}\Big) \right)$, then \begin{equation*}
- {\frac\alpha4}{\frac{\kappa\log2}4} \nu 2^{2\nu} \leq -\alpha\log2\Phi_\lambda'\left( \frac{t+\tau}{\beta} \right)\nu.
\end{equation*} 
On the contrary, if $\nu < {1\over 2\log 2}\log\left(\frac{4\log2}{\kappa'}\Phi_\lambda'\Big( \frac{t+\tau}{\beta}\Big) \right)$ then 
$${
\frac{4\log2}{\kappa'}}\Phi_\lambda'\left( \frac{t+\tau}{\beta} \right) > 2^{2\nu}
$$ 
and, hence, by \eqref{eq:DefPhi}, we obtain \begin{equation*} 
\begin{aligned}
&\frac{1}{4}\frac{t+\tau}{\beta}\Phi_\lambda''\left(\frac{t+\tau}{\beta}\right) = -\frac{1}{4}\lambda \left( \Phi_\lambda'\left(\frac{t+\tau}{\beta}\right) \right)^2 \mu\left( \frac{1}{\Phi_\lambda'\left(\frac{t+\tau}{\beta}\right)} \right) \\
&\qquad \leq -\frac{1}{4}\lambda \left( \Phi_\lambda'\left(\frac{t+\tau}{\beta}\right) \right)^2 \mu\left( \frac{1}{\frac{4\log2}{\kappa'}\Phi_\lambda'\left(\frac{t+\tau}{\beta}\right)} \right) \\
&\qquad\leq -\frac{1}{4}\lambda \frac{\kappa'}{4\log2} \Phi_\lambda'\left(\frac{t+\tau}{\beta}\right)\left( 1+ \log\left( \frac{4\log2}{\kappa'}\Phi_\lambda'\left(\frac{t+\tau}{\beta}\right) \right)  \right) \\
&\qquad\leq -\frac{1}{4}\lambda \frac{\kappa'}{4\log2} \Phi_\lambda'\left(\frac{t+\tau}{\beta}\right) (1+2\nu\log2) \\
&\qquad\leq -\lambda \frac{\kappa'(1+\log2)}{16\log2}  \Phi_\lambda'\left(\frac{t+\tau}{\beta}\right) \nu,
\end{aligned} \end{equation*} 
where we have used the fact  that the function $\varepsilon\mapsto \varepsilon(|\log \varepsilon|+1)$ is increasing. Consequently, if we choose 
$\lambda\geq\bar\lambda$ with 
\begin{equation*}
\bar\lambda \geq \frac{16\alpha (\log2)^2(\sigma+\tau)}{\kappa'(1+\log2)},
\end{equation*} we have 
\begin{equation*}
\frac{1}{4}\frac{t+\tau}{\beta}\Phi''_\lambda\left( \frac{t+\tau}{\beta} \right) \leq -\alpha\log2(t+\tau)\Phi'_\lambda\left( \frac{t+\tau}{\beta} \right)\nu
\end{equation*} and hence, the term \eqref{term1} is compensated by \eqref{term2} and \eqref{term3}.

Now we consider the term 
\begin{equation} \label{term5}
(t+\tau) \alpha^2\log^2(2) \nu^2 \|v_\nu(t)\|_{L^2}.
\end{equation} If $\nu \geq\frac1{\log2} \log\left( \frac{16\alpha\log2}{\kappa} \right) =: \bar{\nu}_1$,  then 
\begin{equation*}
-{\frac\alpha4}\frac{\kappa\log2}{4} \nu 2^{2\nu} + \alpha^2 \log^2(2) \nu^2 \leq 0.
\end{equation*} 
If $\nu \leq \bar{\nu}_1$, then we choose a possibly larger $\bar{\gamma}$ such that 
\begin{equation*}
\frac{\gamma}{4} \geq \alpha^2 \log^2(2) \bar{\nu}_1^2 ( \sigma + \tau )
\end{equation*} 
for all $\gamma \geq \bar{\gamma}$. We obtain 
\begin{equation*}
-\frac{\gamma}{4} + \alpha^2 \log^2(2) \nu^2(t+\tau) \leq 0, 
\end{equation*} 
and, consequently, \eqref{term5} is absorbed by \begin{equation*}
-\frac{\alpha}{4} \frac{\kappa\log2}{4} (t+\tau) \nu 2^{2\nu} \|v_\nu(t)\|_\Lt^2 - \frac{\gamma}{4}\|v_\nu(t)\|_\Lto^2.
\end{equation*} 
The term $- \alpha\gamma \log2 (t+\tau) \nu \|v_\nu(t)\|_\Lto^2$ can be neglected since it is negative. However, we stress here that it is a crucial term in order to achieve our energy estimate for an equation including also lower order terms. Recalling also Propositions \ref{PropBernstein} and \ref{prop:Pos}, we obtain 
\begin{equation*} 
\begin{aligned}
& \frac{1}{2}(t+\tau)\|\partial_t v_\nu(t)\|_{L^2}^2 + \frac{\gamma}{8}\|v_\nu(t)\|_\Lto^2 \\
& \qquad \leq \frac{\gamma}{2}\frac{d}{dt}\left( (t+\tau)\|v_\nu(t)\|_\Lto^2 \right) - \frac{1}{2}\frac{d}{dt}\left( (t+\tau)\Phi_\lambda'\left( \frac{t+\tau}{\beta} \right)\|v_\nu(t)\|_\Lto^2 \right) \\
& \qquad \quad  + \frac{1}{2}\frac{d}{dt}\left( (t+\tau)\la T_{a_\nu}^m \partial_x v_\nu(t) |  \partial_x v_\nu(t) \ral \right) \\
& \qquad \quad - \frac\kappa82^{2\nu}\|v_\nu(t)\|_{L^2}^2\\
& \qquad \quad - \frac{\alpha}{4}\frac{\kappa\log2}{4}(t+\tau)\nu 2^{2\nu}\|v_\nu(t)\|_\Lt^2 \\
& \qquad \quad + \alpha\log2\nu 2^{-\alpha t\nu}(t+\tau)\la v_\nu(t) \mid \mathcal X_\nu(t) \ral \\
& \qquad \quad - (t+\tau)2^{-\alpha t\nu} \la \mathcal X_\nu(t) \mid \partial_t v_\nu(t) \ral
- \frac{\gamma}{8}\|v_\nu(t)\|_\Lto^2.
\end{aligned}\end{equation*}  
Integrating over $[0,s] \subseteq [0,\sigma]$, we get \begin{equation*} \begin{aligned}
& \frac{\kappa}{8} \int_0^s 2^{2\nu }\|v_\nu(t)\|_\Lt^2 dt + \frac{\gamma}{8} \int_0^s \|v_\nu(t)\|_\Lt^2 dt \\
& \qquad \leq \frac12 \tau\Phi_\lambda'\left( \frac{\tau}{\beta} \right)\|v_\nu(0)\|_\Lto^2 + \left(\frac{\gamma}{2}+C^{(4)} 2^{2\nu}\right)(s+\tau)\|v_\nu(s)\|_\Lt^2 \\
& \qquad \quad - \frac{\alpha}{4}\frac{\kappa\log2}{4} \int_0^s (t+\tau)\nu 2^{2\nu}\|v_\nu(t)\|_\Lt^2 dt - \frac{\gamma}{8} \int_0^s \|v_\nu(t)\|_\Lto^2 dt \\
& \qquad \quad  - \frac{1}{2}\int_0^s (t+\tau)\|\partial_t v_\nu(t)\|_{L^2}^2 dt \\
& \qquad \quad - \int_0^s (t+\tau)2^{-\alpha t\nu} \la \mathcal X_\nu(t) \mid \partial_t v_\nu(t) \ral dt \\
& \qquad \quad + \alpha\log2 \int_0^s \nu 2^{-\alpha t\nu}(t+\tau)\la v_\nu(t) \mid \mathcal X_\nu(t) \ral dt,
\end{aligned}\end{equation*} 
where we have used the estimate 
\begin{eqnarray*}
|\la \partial_x v_\nu(s )\mid T_{a_\nu}^m\partial_x v_\nu(s ) \ral| \leq C^{(4)}2^{2\nu}\|v_\nu(s )\|_\Lt^2.
\end{eqnarray*}

\medskip
%%%%%%%%%%%
{\em - End of the proof}
%%%%%%%%%%%
\smallskip

Now we sum over $\nu$ and we obtain \begin{equation*} \begin{aligned}
& \frac{\kappa}{8} \int_0^s \sum\limits_{\nu \geq 0}2^{2\nu}\|v_\nu(t)\|_\Lt^2 dt + \frac{\gamma}{8} \int_0^s \sum\limits_{\nu \geq 0} \|v_\nu(t)\|_\Lt^2 dt \\
& \qquad \leq \frac12 \tau\Phi_\lambda'\left( \frac{\tau}{\beta} \right)\sum\limits_{\nu \geq 0}\|v_\nu(0)\|_\Lto^2 - \frac{\gamma}{8} \int_0^s \sum\limits_{\nu \geq 0} \|v_\nu(t)\|_\Lto^2 dt \\
& \qquad \quad - \frac{1}{2}\int_0^s (t+\tau) \sum\limits_{\nu \geq 0} \|\partial_t v_\nu(t)\|_{L^2}^2 dt \\
& \qquad \quad + \frac{\gamma}{2}(s+\tau)\sum\limits_{\nu \geq 0}\|v_\nu(s)\|_\Lt^2 + C^{(4)} (s+\tau)\sum\limits_{\nu \geq 0} 2^{2\nu}\|v_\nu(s)\|_\Lt^2 \\
& \qquad \quad - \frac{\alpha}{4}{\frac{\kappa\log2}4}\int_0^s (t+\tau)\sum\limits_{\nu \geq 0}\nu 2^{2\nu}\|v_\nu(t)\|_\Lt^2 dt \\
& \qquad \quad - \int_0^s (t+\tau) \sum\limits_{\nu \geq 0} 2^{-\alpha t\nu} \la \mathcal X_\nu(t) \mid \partial_t v_\nu(t) \ral dt \\
& \qquad \quad + \alpha\log2 \int_0^s (t+\tau) \sum\limits_{\nu \geq 0} \nu 2^{-\alpha t\nu}\la v_\nu(t) \mid \mathcal X_\nu(t) \ral dt.
\end{aligned}\end{equation*} 

Now we have
\begin{multline*}
|\sum\limits_{\nu \geq 0} 2^{-\alpha t\nu} \la \mathcal X_\nu(t) \mid \partial_t v_\nu(t) \ral|\\
\leq \sum\limits_{\nu \geq 0} 2^{-\alpha t\nu} |\la \partial_x([\Delta_\nu,T^m_a]\partial_xw(t)) \mid \partial_t v_\nu(t)\ral|\\
+\sum\limits_{\nu \geq 0} 2^{-\alpha t\nu} |\la \Delta_\nu(\partial_x((a-T^m_a)\partial_xw(t))) \mid \partial_t v_\nu(t)\ral|\\
\leq \sum\limits_{\nu \geq 0} 2^{-\alpha t\nu} \| \partial_x([\Delta_\nu,T^m_a]\partial_xw(t))\|_{L^2} \|\partial_t v_\nu(t)\|_{L^2}\\
+\sum\limits_{\nu \geq 0} 2^{-\alpha t\nu} \| \Delta_\nu(\partial_x((a-T^m_a)\partial_xw(t)))\|_{L^2} \|\partial_t v_\nu(t)\|_{L^2}\\
\leq \Big(\sum\limits_{\nu \geq 0} 2^{-2\alpha t\nu} \| \partial_x([\Delta_\nu,T^m_a]\partial_xw(t))\|_{L^2}^2\Big)^{\frac12} \Big(\sum\limits_{\nu \geq 0}\|\partial_t v_\nu(t)\|_{L^2}^2\Big)^{\frac12}\\
+\Big(\sum\limits_{\nu \geq 0} 2^{2(1-\alpha t)\nu} \| \Delta_\nu((a-T^m_a)\partial_xw(t))\|_{L^2}^2\Big)^{\frac12} \Big(\sum\limits_{\nu \geq 0}\|\partial_t v_\nu(t)\|_{L^2}^2\Big)^{\frac12}\\
\end{multline*}
By Corollary \ref{interpolation}, Proposition \ref{comm} and Proposition \ref{SobolevLP} we get
\begin{multline*}
|\sum\limits_{\nu \geq 0} 2^{-\alpha t\nu} \la \mathcal X_\nu(t) \mid \partial_t v_\nu(t) \ral|
\leq C^{(5)}\|w(t)\|_{H^{1-\alpha t}} \Big(\sum\limits_{\nu \geq 0}\|\partial_t v_\nu(t)\|_{L^2}^2\Big)^{\frac12}\\
+C^{(6)}\|(a-T^m_a)\partial_x w(t)\|_{H^{1-\alpha t}} \Big(\sum\limits_{\nu \geq 0}\|\partial_t v_\nu(t)\|_{L^2}^2\Big)^{\frac12}\\
\leq C^{(7)}\|w(t)\|_{H^{1-\alpha t}} \Big(\sum\limits_{\nu \geq 0}\|\partial_t v_\nu(t)\|_{L^2}^2\Big)^{\frac12}\\
\leq C^{(8)}\Big(\sum\limits_{\nu \geq 0}2^{2(1-\alpha t)\nu}\|w_\nu(t)\|_{L^2}^2\Big)^{\frac12} \Big(\sum\limits_{\nu \geq 0}\|\partial_t v_\nu(t)\|_{L^2}^2\Big)^{\frac12}\\
\leq C^{(9)}\Big(\sum\limits_{\nu \geq 0}2^{2\nu}\|v_\nu(t)\|_{L^2}^2\Big)^{\frac12} \Big(\sum\limits_{\nu \geq 0}\|\partial_t v_\nu(t)\|_{L^2}^2\Big)^{\frac12}\\
\leq C^{(10)}\sum\limits_{\nu \geq 0}2^{2\nu}\|v_\nu(t)\|_{L^2}^2+\frac12\sum\limits_{\nu \geq 0}\|\partial_t v_\nu(t)\|_{L^2}^2.
\end{multline*}
In the same way one can prove that
\begin{equation*}
|\sum\limits_{\nu \geq 0}\nu 2^{-\alpha t\nu} \la \mathcal X_\nu(t) \mid  v_\nu(t) \ral|\leq  C^{(11)}\sum\limits_{\nu \geq 0}2^{2\nu}\|v_\nu(t)\|_{L^2}^2.
\end{equation*}
We thus obtain
\begin{equation*} \begin{aligned}
& \frac{\kappa}{8} \int_0^s \sum\limits_{\nu \geq 0}2^{2\nu}\|v_\nu(t)\|_\Lt^2 dt + \frac{\gamma}{8} \int_0^s \sum\limits_{\nu \geq 0} \|v_\nu(t)\|_\Lt^2 dt \\
& \qquad \leq \frac12 \tau\Phi_\lambda'\left( \frac{\tau}{\beta} \right)\sum\limits_{\nu \geq 0}\|v_\nu(0)\|_\Lto^2 - \frac{\gamma}{8} \int_0^s \sum\limits_{\nu \geq 0} \|v_\nu(t)\|_\Lto^2 dt \\
& \qquad \quad + \frac{\gamma}{2}(s+\tau)\sum\limits_{\nu \geq 0}\|v_\nu(s)\|_\Lt^2 + C^{(4)} (s+\tau)\sum\limits_{\nu \geq 0} 2^{2\nu}\|v_\nu(s)\|_\Lt^2 \\
& \qquad \quad - \frac{\alpha}{4}{\frac{\kappa\log2}4}\int_0^s (t+\tau)\sum\limits_{\nu \geq 0}\nu 2^{2\nu}\|v_\nu(t)\|_\Lt^2 dt \\
& \qquad \quad + C^{(12)} \int_0^s (t+\tau) \sum\limits_{\nu \geq 0} 2^{2\nu}\|v_\nu(t)\|_{L^2}^2  dt.
\end{aligned}\end{equation*} 
Now the term 
$$
C^{(12)} \int_0^s (t+\tau) \sum\limits_{\nu \geq 0} 2^{2\nu}\|v_\nu(t)\|_{L^2}^2  dt
$$
can be absorbed by 
$$
- \frac{\alpha}{4}{\frac{\kappa\log2}4}\int_0^s (t+\tau)\sum\limits_{\nu \geq 0}\nu 2^{2\nu}\|v_\nu(t)\|_\Lt^2 dt
$$
for high frequencies, and by 
$$
- \frac{\gamma}{8} \int_0^s \sum\limits_{\nu \geq 0} \|v_\nu(t)\|_\Lto^2 dt 
$$
for low frequencies by choosing $\bar\gamma$ larger if necessary.

All in all, we finally obtain
\begin{equation*} 
\begin{aligned}
& \frac{\kappa}{8} \int_0^s \sum\limits_{\nu \geq 0} 2^{2\nu}\|v_\nu(t)\|_\Lt^2 dt + \frac{\gamma}{8} \int_0^s \sum\limits_{\nu \geq 0} \|v_\nu(t)\|_\Lt^2 dt \\
& \qquad \leq \frac12 \tau\Phi_\lambda'\left( \frac{\tau}{\beta} \right)\sum\limits_{\nu \geq 0}\|v_\nu(0)\|_\Lto^2 + \frac{\gamma}{2}(s+\tau)\sum\limits_{\nu \geq 0}\|v_\nu(s)\|_\Lt^2 \\
& \qquad \quad + C^{(4)} (s+\tau)\sum\limits_{\nu \geq 0} 2^{2\nu}\|v_\nu(s)\|_\Lt^2.
\end{aligned}
\end{equation*} 
From this, going back to $u_\nu$ and using Proposition \ref{SobolevLP}, the weighted energy estimate \eqref{EnergyEst1} follows. \qed

%%%%%%%%%%%%%%%%%%%%%%%%%%%%%%%%%%%%%%%%%%%%%%%%%%%%%%%%%%%%%%%%
\section{Conditional stability up to the final time} 
%%%%%%%%%%%%%%%%%%%%%%%%%%%%%%%%%%%%%%%%%%%%%%%%%%%%%%%%%%%%%%%%
In this section we state and prove two global stability theorems for solutions of (\ref{BWeq}) up to the final time $T$. The first result gives a logarithmic type control of $\|u\|_{L^2((0,T),L^2)}$ in terms of $\|u(0)\|_{L^2}$.

\begin{thm}\label{L2final}
Assume Hypothesis \ref{mainhypothesis} is satisfied.
Then for all $D_0>0$ there exist positive constants $\rho''$, $\delta''$ and $K''$, depending only on $A_{LL}$, $A$, $\kappa$, $T$ and $D_0$, such that if $u \in \mathcal H$ is a solution of \eqref{BWeq} satisfying $\sup_{t\in [0,T]}\|u(t,\cdot)\|_{L^2} \leq D_0$ and $\|u(0,\cdot)\|_{L^2} \leq \rho''$, the inequality \begin{eqnarray*}
\|u\|_{L^2((0,T),L^2)} \leq K'' \frac1{|\log\|u(0)\|_{L^2}|^{\delta''}}
\end{eqnarray*} holds true.\end{thm}
\begin{rem}Notice that, following Remark \ref{forward}, it would be sufficient to impose an a-priory bound on $\|u(T,\cdot)\|_{L^2}$, which authomatically implies the a-priori bound for $\|u(t,\cdot)\|_{L^2}$, $t\in[0,T]$. \end{rem}

\begin{proof}[Proof of Theorem {\ref{L2final}}]
First we observe that, due to Theorem \ref{thm:DSP2009global}, it is not restrictive to assume that $\alpha_1\leq T^{-1}$. Indeed, if this is not the case we can take $T'$, $0<T'<T$, such that $T-T'<\alpha_1^{-1}$, and then in $[0,T']$ we apply the pointwise estimate given by Theorem \ref{thm:DSP2009global}, so we just need to estimate $\int_{T'}^T\|u(t)\|^2_{L^2}\,dt$ in terms of $\|u(T')\|_{L^2}$. 

With such assumption we can apply Proposition \ref{thm:Energy} with $\alpha=1/T$, $\sigma=T$ and $\tau=T/4$ and we can find $\lambda>1$, $\gamma>0$ and $M>0$ such that for all $\beta\geq T+\tau=\frac54\tau$ and whenever $u \in \mathcal H$ is a solution of equation \eqref{BWeq}, then

\begin{multline*}
\int_0^T e^{2\gamma t}e^{-2\beta \Phi_\lambda({t+\tau\over
\beta})} \|u(t,\cdot)\|^2_{H^{1-\alpha t}}\, dt\\ \leq
M\gamma((T+\tau) e^{2\gamma T}e^{-2\beta
\Phi_\lambda({T+\tau\over \beta})}\|u(T,\cdot)\|^2_{L^2}
+\tau\,\Phi_\lambda'({\tau\over\beta})e^{-2\beta
\Phi_\lambda({\tau\over \beta})}\|u(0,\cdot)\|^2_{L^2}).
\end{multline*}
Now for any $r\in (0,T)$ we have
\begin{multline*}
\int_0^{T-r} e^{2\gamma t}e^{-2\beta \Phi_\lambda({t+\tau\over
\beta})} \|u(t,\cdot)\|^2_{L^2}\, dt\\ \leq
M\gamma((T+\tau) e^{2\gamma T}e^{-2\beta
\Phi_\lambda({T+\tau\over \beta})}\|u(T,\cdot)\|^2_{L^2}
+\tau\,\Phi_\lambda'({\tau\over\beta})e^{-2\beta
\Phi_\lambda({\tau\over \beta})}\|u(0,\cdot)\|^2_{L^2}),
\end{multline*}
where we have used the fact that $\|u(t,\cdot)\|_{L^2}\leq
\|u(t,\cdot)\|_{H^{1-\alpha t}}$. Now, the
function $\Phi_\lambda$ is increasing and consequently the
function $t\mapsto e^{-2\beta \Phi_\lambda({(t+\tau)/ \beta})}$ is
decreasing. We deduce that
\begin{multline*}
 e^{-2\beta \Phi_\lambda({T-r+\tau\over \beta})}
\int_0^{T-r}\|u(t,\cdot)\|^2_{L^2}\,dt\\ \leq M'\left(e^{-2\beta \Phi_\lambda({T+\tau\over
\beta})}\|u(T,\cdot)\|^2_{L^2}
+\Phi_\lambda'({\tau\over\beta})e^{-2\beta
\Phi_\lambda({\tau\over \beta})}\|u(0,\cdot)\|^2_{L^2}\right),
\end{multline*}
where $M'=M\gamma2Te^{2\gamma T}$.
Then
\begin{multline*}
\int_0^{T-r}\|u(t,\cdot)\|^2_{L^2}\,dt \leq  M'
\Phi_\lambda'({\tau\over\beta})\Big(e^{2\beta(
\Phi_\lambda({T-r+\tau\over
\beta})-\Phi_\lambda({T+\tau\over
\beta}))}\|u(T,\cdot)\|^2_{L^2}
\\
 +e^{2\beta( \Phi_\lambda({T-r+\tau\over
\beta})-\Phi_\lambda({\tau\over
\beta}))}\|u(0,\cdot)\|^2_{L^2}\Big)
\\
\leq   M'
\Phi_\lambda'({\tau\over\beta})e^{2\beta(
\Phi_\lambda({T-r+\tau\over
\beta})-\Phi_\lambda({T+\tau\over
\beta}))}\Big(\|u(T,\cdot)\|^2_{L^2}
\\
 +e^{-2\beta \Phi_\lambda({\tau\over
\beta})}\|u(0,\cdot)\|^2_{L^2}\Big),
\end{multline*}
where we used the fact that $\Phi_\lambda'(\frac\tau\beta)\geq 1$ and $\Phi_\lambda(\frac{T+\tau}\beta)\leq0$.
We recall that the function $\Phi_\lambda$ is concave, so 
\begin{multline*}
\Phi_\lambda({T-r+\tau\over \beta})-\Phi_\lambda({T+\tau\over
\beta})\\
\leq \Phi'_\lambda({T+\tau\over \beta})({T-r+\tau\over \beta}-{T+\tau\over \beta})=-\Phi'_\lambda({T+\tau\over \beta}){r\over \beta},
\end{multline*}
and then
\begin{multline*}
\int_0^{T-r}\|u(t,\cdot)\|^2_{L^2}\,dt \\
\leq   M'
\Phi_\lambda'({\tau\over\beta})e^{-2r\Phi'_\lambda({T+\tau\over \beta})}\Big(\|u(T,\cdot)\|^2_{L^2}
 +e^{-2\beta \Phi_\lambda({\tau\over
\beta})}\|u(0,\cdot)\|^2_{L^2}\Big).
\end{multline*}
By Lemma \ref{weight-properties} we have that
$$
\Phi'_\lambda({T+\tau\over \beta})=
\psi_\lambda({T+\tau\over \tau}{\tau\over \beta})=
\exp\Big({\Big({T+\tau\over \tau}\Big)^{-\lambda}-1}\Big)\Big(\psi_\lambda({\tau\over \beta})\Big)^{({T+\tau\over \tau})^{-\lambda}}.
$$
We remind that $\tau=T/4$, so $\frac{T+\tau}\tau=5$, and 
$$
\Phi'_\lambda({T+\tau\over \beta})=\bar N\psi_\lambda(\frac\tau\beta)^{\bar\delta},
$$
where $\bar\delta=5^{-\lambda}$ and $\bar N=e^{\bar\delta-1}$.
It follows that 
\begin{multline*}
\int_0^{T-r}\|u(t,\cdot)\|^2_{L^2}\,dt \\
\leq   M'
\psi_\lambda({\tau\over\beta})e^{-2r\bar N\psi_\lambda(\frac\tau\beta)^{\bar\delta}}\Big(\|u(T,\cdot)\|^2_{L^2}
 +e^{-2\beta \Phi_\lambda({\tau\over
\beta})}\|u(0,\cdot)\|^2_{L^2}\Big).
\end{multline*}

Now we observe that 
$$
\psi_\lambda({\tau\over\beta})e^{-r\bar N\psi_\lambda(\frac\tau\beta)^{\bar\delta}}=
r^{-1/\bar\delta}r^{1/\bar\delta}\psi_\lambda({\tau\over\beta})e^{-r\bar N\psi_\lambda(\frac\tau\beta)^{\bar\delta}}
\leq C_{\bar N,\bar\delta}\, r^{-1/\bar\delta}
$$
where 
$$
C_{\bar N,\bar\delta}:=\sup_{z\geq 0} ze^{-\bar N z^{\bar\delta}}.
$$
Then 
\begin{multline*}
\int_0^{T-r}\|u(t,\cdot)\|^2_{L^2}\,dt \\
\leq   M'C_{\bar N,\bar\delta}\, r^{-1/\bar\delta}
e^{-r\bar N\psi_\lambda(\frac\tau\beta)^{\bar\delta}}\Big(\|u(T,\cdot)\|^2_{L^2}
 +e^{-2\beta \Phi_\lambda({\tau\over
\beta})}\|u(0,\cdot)\|^2_{L^2}\Big).
\end{multline*}

We choose now $\beta$ in such a way that $e^{-\beta \Phi_\lambda({\tau\over \beta})}=
\|u(0,\cdot)\|^{-1}_{L^2}$ i. e.
$$
{\beta\over \tau} \Phi_\lambda({\tau\over \beta})={1\over \tau}\log \|u(0,\cdot)\|_{L^2}.
$$
 We obtain
$ \beta= \tau \Lambda_\lambda^{-1}({1\over \tau}\log \|u(0,\cdot)\|_{L^2})$, where $\Lambda_\lambda(y)=y\Phi_\lambda(1/y)$. If  $\|u(0,\cdot)\|_{L^2}\leq \bar \rho:= e^{\tau\Lambda_\lambda(5)}$, then $\beta\geq T+\tau$. We have then
\begin{multline*}
\int_0^{T-r}\|u(t,\cdot)\|^2_{L^2}\,dt \\
\leq   M'C_{\bar N,\bar\delta}\, r^{-1/\bar\delta}
e^{-r\bar N\psi_\lambda(\frac1{\Lambda_\lambda^{-1}({1\over \tau}\log \|u(0,\cdot)\|_{L^2})})^{\bar\delta}}\Big(\|u(T,\cdot)\|^2_{L^2}
 +1\Big).
\end{multline*}
By Lemma \ref{weight-properties} we have that
\begin{equation*}
\lim_{z \rightarrow -\infty} -\frac{1}{z}\psi_\lambda\big( \frac{1}{\Lambda_\lambda^{-1}(z)} \big) = +\infty,\end{equation*}
so 
$$
\psi_\lambda\big( \frac{1}{\Lambda_\lambda^{-1}(z)} \big) \geq|z|
$$
if $z<0$ and $|z|$ is sufficiently large. It follows that there exists $\tilde\rho\leq\bar\rho$ such that, if $\|u(0)\|_{L^2}\leq\tilde\rho$, then
\begin{multline*}
\int_0^{T-r}\|u(t,\cdot)\|^2_{L^2}\,dt \\
\leq   M'C_{\bar N,\bar\delta}\, r^{-1/\bar\delta}
e^{-r\bar N(\frac1\tau |\log \|u(0,\cdot)\|_{L^2})|)^{\bar\delta}}
\Big(\|u(T,\cdot)\|^2_{L^2}
 +1\Big).
\end{multline*}
On the other hand,
\begin{equation*}
\int_{T-r}^T\|u(t,\cdot)\|^2_{L^2}\,dt \leq D_0 r.
\end{equation*}
It follows that for all $r>0$
\begin{equation*}
\int_0^{T}\|u(t,\cdot)\|^2_{L^2}\,dt 
\leq   \bar M\big(r+ r^{-1/\bar\delta}
e^{-r\tilde N( |\log \|u(0,\cdot)\|_{L^2})|)^{\bar\delta}}
\Big),
\end{equation*}
where
$$
\bar M=(M'C_{\bar N,\bar\delta}+1)(D_0+1) \quad\text{and}\quad \tilde N=\frac{\bar N}{\tau^{\bar\delta}}.
$$
Finally we choose 
$$
r=|\log\|u(0)\|_{L^2}|^{-\bar\delta/2},
$$
so we get
\begin{multline*}
\int_0^{T}\|u(t,\cdot)\|^2_{L^2}\,dt \\
\leq   \bar M\big(\frac1{|\log\|u(0)\|_{L^2}|^{\bar\delta/2}}+
|\log\|u(0)\|_{L^2}|^{1/2}
e^{-\tilde N( |\log \|u(0,\cdot)\|_{L^2})|)^{\bar\delta/2}}
\Big)\\
\leq\bar M(1+E_{\tilde N,\bar\delta})\frac1{|\log\|u(0)\|_{L^2}|^{\bar\delta/2}},
\end{multline*}
where
$$
E_{\tilde N,\bar\delta}=\sup_{z\geq 0}z^{(1+\bar\delta)/2}e^{-\tilde Nz^{\bar\delta/2}}.
$$
The proof is complete.
\end{proof}

Under a stronger a-priori bound on admissible solutions in $[0,T]$, namely assuming an a-priori bound in $H^1$ rather than in $L^2$, we can prove a pointwise stability estimate of logarithmic type up to the final time $T$.

\begin{thm}\label{Linftyfinal}
Assume Hypothesis \ref{mainhypothesis} is satisfied.
Then for all $D_1>0$ there exist positive constants $\rho'''$, $\delta'''$ and $K'''$, depending only on $A_{LL}$, $A$, $\kappa$, $T$ and $D_1$, such that if $u \in \mathcal H$ is a solution of \eqref{BWeq} satisfying $\sup_{t\in [0,T]}\|u(t,\cdot)\|_{H^1} \leq D_1$ and $\|u(0,\cdot)\|_{L^2} \leq \rho'''$, the inequality \begin{eqnarray*}
\sup_{t\in[0,T]}\|u(t,\cdot)\|_{L^2} \leq K''' \frac1{|\log\|u(0,\cdot)\|_{L^2}|^{\delta'''}}
\end{eqnarray*} holds true. \end{thm}
\begin{rem}Notice that, following Remark \ref{forward}, it would be sufficient to impose an a-priory bound on $\|u(T,\cdot)\|_{H^1}$, which authomatically implies the a-priori bound for $\|u(t,\cdot)\|_{H^1}$, $t\in[0,T]$. \end{rem}
\begin{proof}[Proof of Theorem {\ref{Linftyfinal}}]
We begin by noticing that, since $u$ solves (\ref{BWeq}), then
$$
\|\partial_t u(t,\cdot)\|_{H^{-1}}\leq\frac1\kappa D_1
$$
It follows from Morrey's inequality that
$$
\sup_{t\in[0,T]}\|u(t,\cdot)\|_{H^{-1}}\leq C_T\|u\|_{L^2((0,T), H^{-1})}^{1/2}\|u\|_{H^1((0,T), H^{-1})}^{1/2}
$$
(for a direct simple proof see \cite[proof of Thm. 8.8]{brezis}).
Then by Theorem \ref{L2final} for $\|u(0)\|_{L^2}\leq\rho''$ we get
$$
\sup_{t\in[0,T]}\|u(t,\cdot)\|_{H^{-1}}\leq C_T\left(\frac1\kappa D_1\right)^{1/2}\left( K'' \frac1{|\log\|u(0)\|_{L^2}|^{\delta''}}\right)^{1/2}.
$$
The conclusion follows observing that for each fixed $t\in[0,T]$ we have
\begin{multline*}
\|u(t,\cdot)\|_{L^2}\leq\|u(t,\cdot)\|_{H^1}^{1/2}\|u(t,\cdot)\|_{H^{-1}}^{1/2}\\
\leq D_1^{1/2}C_T^{1/2}\left(\frac1\kappa D_1\right)^{1/4}\left( K'' \frac1{|\log\|u(0)\|_{L^2}|^{\delta''}}\right)^{1/4}.
\end{multline*}
\end{proof}

%%%%%%%%%%%%%%%%%%%%%%%%%%%%%%%%%%%%%%%%%%%%%%%%%%%%%%%%%%%%%%%%
\section{Reconstruction of the initial condition for parabolic equations} 
%%%%%%%%%%%%%%%%%%%%%%%%%%%%%%%%%%%%%%%%%%%%%%%%%%%%%%%%%%%%%%%%
In view of applications it is convenient to rephrase Theorem \ref{Linftyfinal}. Consider the (forward) parabolic equation
\begin{equation}
\label{FWeq}
 \partial_t u - \sum\limits_{j,k=1}^n \partial_{x_j}(a_{jk}(t,x)\partial_{x_k}u) = 0
\end{equation} on the strip $[0,T] \times \R^n_x$ and assume Hypothesis \ref{mainhypothesis} is satisfied. Then we have:

\begin{cor}\label{reconst}
Let $D>0$. There exist positive constants $\rho_D$, $\delta_D$ and $K_D$, depending only on $A_{LL}$, $A$, $\kappa$, $T$ and $D$, such that if $u,v \in C^0([0,T],H^1)\cap C^1(]0,T],L^2)$ are solutions of \eqref{FWeq} satisfying $\|u(0,\cdot)\|_{H^1}\leq D$, $\|v(0,\cdot)\|_{H^1}\leq D$ and $\|u(T,\cdot)-v(T,\cdot))\|_{L^2} \leq \rho_D$, then the inequality \begin{eqnarray*}
\|u(0,\cdot)- v(0,\cdot)\|_{L^2} \leq K_D \frac1{|\log\|u(T,\cdot)-v(T,\cdot)\|_{L^2}|^{\delta_D}}
\end{eqnarray*} holds true. \end{cor}

Corollary \ref{reconst} can be exploited to reconstruct the initial condition of an unknown solution $u(t)$ of (\ref{FWeq}), provided we can measure with arbitrary accuracy its final configuration $u_T:=u(T)$. More precisely, suppose that for every $\theta>0$ we can perform a measurement $v_{\theta,T}$ of $u_T$ such that 
\begin{equation*}
\|v_{\theta,T}-u_T\|_{L^2}\leq\theta.
\end{equation*}
Moreover, suppose that we know {\em a priori} that $\|u(0)\|_{H^1}\leq D$ for some $D>0$. We are interested in finding a computable approximation of $u(0)$. If it were possible to solve equation (\ref{FWeq}) backward in time with final condition $v(T)=v_{\theta,T}$, then by Corollary \ref{reconst} we would get that $v(0)$ is closed to $u(0)$, provided $\|v(0)\|_{H^1}\leq D$ and $v_{\theta,T}$ is sufficiently closed to $u_T$. However, equation (\ref{FWeq}) with final condition $v(T)=v_{\theta,T}$ in general has no solution, due to the regularizing effect of equation (\ref{FWeq}) forward in time, and to the fact that $v_{\theta,T}$ does not possess any regularity, since it is the output of a measurement. There are various strategies to overcome this major obstruction. We mention the technique of quasi reversibility (see e.g. \cite{Ew}), which consists in perturbing the equation to make it solvable backward in time, and the technique of Fourier truncation, which consists in approximating $v_{\theta,T}$ with a very regular function obtained truncating its Fourier transform. We illustrate the second technique through an example inspired by \cite{FXQ} (see also \cite{Hao}). 

We consider the equation
\begin{equation}\label{simpleFWeq}
 \partial_t u - \sum\limits_{j,k=1}^n a_{jk}(t)\partial_{x_j}\partial_{x_k}u = 0
\end{equation} on the strip $[0,T] \times \R^n_x$ and assume that the coefficients $a_{jk}(t)$ are Log-Lipschitz continuous. Moreover, setting
$$
a(t,\xi):=\sum\limits_{j,k=1}^n a_{jk}(t)\xi_j\xi_k,
$$
we assume that 
$$
\frac12|\xi|^2\leq a(t,\xi)\leq2|\xi|^2,\quad (t,\xi)\in[0,T]\times \R^n_\xi.
$$
Denote by $\mathcal F$ the Fourier transform with respect to the $x$ variable, and by $\mathcal F^{-1}$ its inverse. Let $u\in C^0([0,T],H^1)\cap C^1(]0,T],L^2)$ be a solution of (\ref{simpleFWeq}) and let $\hat u(t,\xi):=(\mathcal F u)(t,\xi)$. Then
$$
\partial_t\hat u(t,\xi)=-a(t,\xi)\hat u(t,\xi).
$$
We set
$$
A(t,\xi):=\int_0^ta(s,\xi)\,ds,
$$
and we observe that $A(t,\xi)$ is increasing in $t$. Since $u(0,\cdot)\in L^2(\R^n_x)$, we have the following explicit representation of $\hat u(t,\xi)$ and hence of $u(t,x)$:
$$
\hat u(t,\xi)=e^{-A(t,\xi)}\hat u(0,\xi), \quad (t,\xi)\in [0,T]\times\R^n_\xi.
$$
 On the othe hand, if $\phi_T(\xi)$ is such that 
 \begin{equation}\label{existcond}e^{A(T,\xi)}\phi_T(\xi)\in L^2(\R^n_\xi), \end{equation}
 then (\ref{simpleFWeq}) can be solved backward in time with final condition $w(T)=w_T:=\mathcal F^{-1}\phi_T$ and the explicit solution is $w(t,x)= (\mathcal F^{-1}\phi)(t,x)$, where
 $$
 \phi(t,\xi)=e^{A(T,\xi)-A(t,\xi)}\phi_T(\xi).
 $$
As above, suppose we know {\em a priori} that $\|u(0)\|_{H^1}\leq D$. Moreover,  suppose that for every $\theta>0$ we can perform a measurement $v_{\theta,T}$ of $u_T$ such that 
\begin{equation*}
\|v_{\theta,T}-u_T\|_{L^2}\leq\theta.
\end{equation*}
 Let $\hat u_T$ and $\hat v_{\theta,T}$ be the Fourier transform of $u_T$ and $v_{\theta,T}$. For $R>0$ define:
$$
\hat u_{T,R}(\xi):=\chi_R(\xi)\hat u_T(\xi) \quad \text{and} \quad \hat v_{\theta,T,R}(\xi):=\chi_R(\xi)\hat v_{\theta,T}(\xi)
$$
 where $\chi_R(\xi)$ is the characteristic function of the ball of radius $R$ in $\>R^n_\xi$. Both $\hat u_{T,R}$ and $\hat v_{\theta,T,R}$
 satisfy (\ref{existcond}) so we can solve (\ref{simpleFWeq}) backward in time with data at $T$ given by $u_{T,R}=\mathcal F^{-1}\hat u_{T,R}$ and $v_{\theta,T,R}=\mathcal F^{-1}\hat v_{\theta,T,R}$. The explicit representations of the corresponding solutions are $u_R(t,x):=\mathcal F^{-1}(\hat u_{R})(t,x)$ and $v_{\theta,R}(t,x):=\mathcal F^{-1}(\hat v_{\theta,R})(t,x)$, where
\begin{equation*}
\hat u_{R}(t,\xi)=e^{A(T,\xi)-A(t,\xi)}\hat u_{T,R}(\xi)\end{equation*}
and
\begin{equation*}\hat v_{\theta,R}(t,\xi)=e^{A(T,\xi)-A(t,\xi)}\hat v_{\theta,T,R}(\xi).
\end{equation*} 
 It is straightforward to check that $\|u_R(0)\|_{H_1}\leq D$.
 Now we have
\begin{multline*}
\|v_{\theta,R}(0)\|_{H^1}\leq \|u_R(0)\|_{H^1}+ \|v_{\theta,R}(0)-u_R(0)\|_{H^1}\\
 \leq D+\left(\int_{|\xi|\leq R}(1+|\xi|^2)e^{2A(T,\xi)}
|\hat v_{\theta,T,R}(\xi)-\hat u_{T,R}(\xi)|^2\,d\xi\right)^{\frac12}\\
\leq D+(1+R^2)^{\frac12}e^{2TR^2}\|v_{\theta,T}-u_T\|_{L^2}
\leq D+e^{(2T+1)R^2}\theta.
 \end{multline*}
Moreover, we have
\begin{multline*}
\|u_T-v_{\theta,T,R}\|_{L^2}\leq \|\hat u_T-\hat u_{T,R}\|_{L^2}+\|\hat u_{T,R}-\hat v_{\theta,T,R}\|_{L^2}\\
\leq \left(\int_{|\xi|\geq R}|\hat u_T(\xi)|^2\,d\xi\right)^{\frac12}+\left(\int_{|\xi|\leq R}|\hat u_{T}(\xi)-\hat v_{\theta,T}(\xi)|^2 \,d\xi\right)^{\frac12}\\
\leq\left(\int_{|\xi|\geq R}(1+|\xi|^2)(1+|\xi|^2)^{-1}e^{-2A(T,\xi)}|\hat u(0,\xi)|^2\,d\xi\right)^{\frac12} +\theta\\
\leq (1+R^2)^{-\frac12}e^{-(T/2)R^2}\|u(0)\|_{H^1}+\theta
\leq e^{-(T/2)R^2}D+\theta
\end{multline*}
 Now, assuming without loss of generality that $\theta<1$, we choose $R(\theta):= (2T+1)^{-1/2}|\log\theta|^{1/2}$ and we notice that $R(\theta)$ tends to $+\infty$ as $\theta\to 0$. With this choice we have
 $$
 \|v_{\theta,R}(0)\|_{H^1}\leq D+1
 $$
 and 
 $$
 \|u_T-v_{\theta,T,R}\|_{L^2}\leq D\theta^{T/(4T+2)}+\theta\leq (D+1)\theta^{T/(4T+2)}.
 $$
 Now let $\rho=\rho_{D+1}$, $K=K_{D+1}$ and $\delta=\delta_{D+1}$ be the constants given by Corollary \ref{reconst}. Then for sufficiently small $\theta$ we have that 
 $$\|u_T-v_{\theta,T,R}\|_{L^2}\leq \rho.$$
 Finally, by Corollary \ref{reconst}, we get
 $$
 \|u(0)-v_{\theta,R(\theta)}(0)\|_{L^2}\leq \tilde K 
 \frac1{|\log\theta|^\delta},
 $$
 where $\tilde K$ can be explicitly expressed in terms of $T$, $D$, $K$ and $\delta$. Therefore $v_{\theta, R(\theta)}(0)$ is the desired approximation of $u(0)$ in $L^2$.

\end{document}